\documentclass[11pt]{article}

\usepackage{amscd,amssymb}
\usepackage{amsthm,amsmath}
\usepackage{graphicx}
\usepackage{enumerate,amssymb,setspace}
\usepackage{amsfonts,amsthm,amscd}
\usepackage{mathtools}
\usepackage{cite}
\usepackage{mathtools}
\usepackage[all]{xy}

\usepackage[latin1]{inputenc}
\usepackage{tikz}
\usetikzlibrary{shapes,arrows}

\usepackage{caption}
\usepackage[]{subcaption}
\usepackage{tikz}
\usetikzlibrary{snakes}

\DeclareCaptionLabelFormat{opening}{}

\newtheorem{theorem}{Theorem}[section]
\newtheorem{lemma}[theorem]{Lemma}
\newtheorem{prop}[theorem]{Proposition}
\newtheorem{corollary}[theorem]{Corollary}
\newtheorem{conjecture}[theorem]{{Conjecture}}

\newtheorem{definition}[theorem]{{Definition}}
\newtheorem{claim}[theorem]{{Claim}}

\def\bclaim{\begin{claim}}
\def\eclaim{\end{claim}}
\def\bdefin{\begin{definition}}
\def\edefin{\end{definition}}
\def\bcor{\begin{corollary}}
\def\ecor{\end{corollary}}
\def\bthm{\begin{theorem}}
\def\ethm{\end{theorem}}
\def\bconj{\begin{conjecture}}
\def\econj{\end{conjecture}}
\def\blem{\begin{lemma}}
\def\elem{\end{lemma}}
\def\blemma{\begin{lemma}}
\def\elemma{\end{lemma}}
\def\bprop{\begin{prop}}
\def\eprop{\end{prop}}
\def\bremark{\begin{remark}}
\def\eremark{\end{remark}}

\theoremstyle{remark}
\newtheorem{remark}[theorem]{Remark}

\def\PR#1{Proposition \ref{#1}}

\def\med{\medskip}

\newcommand{\Pic}{\mathrm{Pic}}

\newcommand{\rk}{\mathrm{rk}}

\newcommand{\mult}{\mathrm{mult}}

\makeatletter\@addtoreset{equation}{section} \makeatother

\font\sml=cmr6  
\newcommand{\dcal}{\mathcal{D}}
\newcommand{\ecal}{\mathcal{E}}

\newcommand{\ocal}{\mathcal{O}}

  \def\calE{\ecal}
 \def\calO{\ocal}
\def\calD{\dcal}

\def\a{\alpha}
\def\al{\alpha}
\def\be{\beta}

 \def\eps{\epsilon}
\def\K{K\"ahler } 
\def\KEE{K\"ahler--Einstein edge } 
\def\h#1{\hbox{#1}}

\def\q{\quad} 
\def\ra{\rightarrow}

\def\Fn{\mathbb{F}_n}
\def\Pic{\operatorname{Pic}}
\def\sm{\setminus}

\newcommand{\PP}{{\mathbb P}} \newcommand{\RR}{\mathbb{R}}
 
 \newcommand{\NN}{{\mathbb N}}
\newcommand{\FF}{{\mathbb F}}
\def\Fn{\FF_n}

\def\beq{\begin{equation}}
\def\eeq{\end{equation}}
\def\bpf{\begin{proof}}
\def\epf{\end{proof}}
\def\bremark{\begin{remark}}
\def\eremark{\end{remark}}

\def\eaeq{\end{aligned}}
\def\baeq{\begin{aligned}}
\def\mult{\operatorname{mult}\!}
\def\saldp{strongly asymptotically log del Pezzo }

\def\saldps{strongly asymptotically log del Pezzos }

\def\aldp{asymptotically log del Pezzo }

\def\aldps{asymptotically log del Pezzos }

\def\saldpno{strongly asymptotically log del Pezzo}

\def\aldpno{asymptotically log del Pezzo}

\def\alf{asymptotically log  Fano }

\def\lb{\label}
\def\er{\eqref}
\def\noi{\noindent}

\newcommand\blfootnote[1]{%
	\begingroup
	\renewcommand\thefootnote{}\footnote{#1}%
	\addtocounter{footnote}{-1}%
	\endgroup
}

\def\Sig{\Sigma}

\def\sign{\h{\rm sign}}

\def\tSig{\tilde \Sigma}

\def\tc{{\tilde c}}

\def\calEprop{\calE_{\perp C}}
\def\calEaway{\calE_{\sm C}}
\def\calEC{\calE_{C}}

\begin{document}
\title{Classification of strongly asymptotically 
\\
log del Pezzo flags and surfaces}
\author{Yanir A. Rubinstein}
\date{}

\maketitle

\begin{abstract}
We introduce the notion of strongly asymptotically log del Pezzo flags,
and classify such flags under the assumption
that their zero-dimensional part lies in the boundary.
We use this result to give a new and conceptual proof of the classification of strongly 
asymptotically log del Pezzo surfaces, originally due to Cheltsov and the author. 
\end{abstract}

\blfootnote{It is a pleasure to thank I.A. Cheltsov for many helpful discussions over the years and our collaboration \cite{CR} which stimulated this article and the referees for
an extremely careful reading and numerous corrections and comments. 
This research was supported by
		NSF grants DMS-1515703,1906370,2204347, a Sloan Research Fellowship, and the Rosi \& Max Varon Visiting
		Professorship in Fall 2019 and Spring 2020 at the Weizmann Institute of Science to which the author is
		grateful, and particularly to B. Klartag and S. Yakovenko, for the excellent research conditions during that time. 
	}
	
\section{Introduction}
\label{}

The classification problem for smooth Fano manifolds of low dimension
has been a fundamental problem in algebraic geometry starting with
the work of the Italian School 
in the second half of the 19th century.  
In an attempt to generalize the problem to {\it pairs}, Maeda introduced
the following notion \cite{Maeda0,Maeda}.

\begin{definition}
\label{definition:log-del-Pezzo-pair} 
Let $X$ be a smooth variety and let $D$ be a simple normal
crossing divisor in $X$.
The
pair $(X,D=\sum D_i)$ is called {log Fano} if
$-K_X-D$ is ample.
\end{definition}
\noindent

Maeda posed the problem of classifying log Fano pairs and gave
a complete classification up to dimension 3
(see also Loginov \cite{Loginov}), which we will come back to shortly. 

Motivated by the study of \KEE metrics, Cheltsov and the author introduced the
following generalization of Maeda's notion, that allows for the coefficients of
the $D_i$ to be slightly less than $1$ \cite[Definition 1.1]{CR}:

\begin{definition}
		\label{definition:log-Fano}
		We say that a pair $(X,D)$ consisting
		of a smooth
	    complex variety $X$ and a divisor $D=\sum_{i=1}^rD_i$ with simple normal crossings
		on $X$ is {\it strongly asymptotically log Fano} if there exists $\eps>0$ such that 
		$-K_X-\sum_{i=1}^r(1-\beta_i)D_i$ is ample for all $(\beta_1,\ldots,\beta_r)\in(0,\eps]^r$.
	\end{definition}

In fact, Maeda's notion corresponds to the case $\be_1=\cdots=\be_r=0$, which by openness of ampleness \cite[Example 1.3.14]{Laz},
implies ampleness for small $\be_i$. In other words, every log Fano is 
strongly asymptotically log Fano. As we will see below the converse is far from true.
For a survey on Maeda's and Cheltov--Rubinstein's work we refer 
to \cite[\S8]{R14}.

The differential geometric interpretation of positive $\be_i$'s is as follows. 
By the resolution of the Calabi--Tian conjecture \cite[Theorem 2]{JMR}, \cite{GP,MR}  
$(X,D)$ is strongly \alf if and only if it admits \K metrics with positive Ricci curvature on $X\sm D$ 
and with edge singularities of angle $2\pi\be_i$
along $D_i$ for all small $\be_i$. Edge singularities are, roughly, conic singularities transverse to the 
`edges' $D_i$. We refer to \cite[\S4]{R14} for a detailled survey.

As customary, from now and on, when discussing dimension 2 we will use `del Pezzo' instead of `Fano', and replace $(X,D)$ by $(S,C)$).
Already in dimension 2, Maeda's classification shows that it is rare for a pair
to be log del Pezzo  \cite[\S3]{Maeda} (see \cite[Proposition 4.1]{CMGR}
for an expository proof):

	\begin{prop}\label{MaedaThm}
		Log del Pezzo pairs $(S,C)$ are classified as follows:
		\hfill\break
		(i) $S=\mathbb{P}^2$, and $C$ is a line,%
		\hfill\break
		(ii)
		$S=\mathbb{P}^2$, and $C=C_1+C_2$, where each $C_i$ is a line,%
		\hfill\break
		(iii)
		$S=\mathbb{P}^2$, and $C$ is a smooth conic,%
		\hfill\break
		(iv)
		$S=\mathbb{F}_{n}$ for some $n\in\NN\cup\{0\}$, and 
		$C$ is a $-n$-curve,
		\hfill\break
		(v)
		$S\cong\mathbb{F}_n$ for any $n\geq 0$, and $C=C_1+C_2$
		where $C_1$ is a $-n$-curve and $C_2$ is a 0-curve,
		\hfill\break
		(vi)
		$S=\mathbb{F}_1$, and $C$ is a smooth $1$-curve,
		\hfill\break
		(vii)
		$S=\mathbb{P}^1\times\mathbb{P}^1$, and $C$ is a smooth 
		$2$-curve.
		
	\end{prop}

In particular, the rank of the Picard group is at most 2 which is extremely restrictive. 
What is more, the notion of log del Pezzo pairs does not even recover the
classical notion of del Pezzo surfaces. 
As it turns out, strongly \aldps are rather richer in structure and the rank of their Picard group has no upper bound. 
From a geometric point of view this is interesting since it provides a host of manifolds on which to construct
canonical \K edge metrics as well as metrics of positive Ricci curvature away from $C$. 
In addition, one easily realizes the classical
del Pezzo surfaces in this picture by considering $(S,C)$ with $C$ a smooth 
anticanonical curve on a del Pezzo surface $S$ as then $-K_S-(1-\be)C\sim -\be K_S>0$.
	
Cheltsov and the author gave a classification of \saldp pairs based on 
two main steps
\cite[Theorems 2.1,3.1]{CR}. First, an induction on
the rank of the Picard group and the key observation that every \saldp pair can be obtained from a pair $(s,c)$
with Picard group at most 2 via blow-ups along smooth points of the boundary $c$ and replacing
the original boundary with its proper transform $\tc$. 
Second, an ad hoc verification of the generality conditions on the blow-ups
allowed for each of the resulting pairs $(s,c)$. 

Here,  we would like to present another
proof of the classification of \saldp pairs, but
our main goal is to actually introduce a different
point of view of independent interest and with
a number of applications. 
The first, necessary, part of the proof is quite
similar to the one in \cite{CR} though we give a more 
conceptual/pedagogical argument (see the flowchart, Figure \ref{Figure0}, in the proof
of Proposition \ref{noawaynotailProp})
that has the advantage of
generalizing to the setting of \aldpno. 
For the second, sufficient, part of the proof,
we present a new approach 
via the notion of SALdP flags that we introduce here
(Definition \ref{badcurvebetadef}).
We show that this notion gives a new characterization of
asymptotically log del Pezzo pairs 
that uses in a precise way blow-downs of the pair.
This characterization is the first main result of this article
(Theorem \ref{badflagcriterionprop}).
The classification of all such flags (Theorem \ref{MainFlagthm}) 
is the second main result of the present note.
We apply this result
to give a slightly cleaner picture of the generality conditions on the blown-up points on the boundary
in the classification of \saldp pairs (Theorem  \ref{mainthm}). 
We believe this notion should also be important for the classification of
the much larger class of \aldp pairs \cite{CR,R14} and to the structure
of the body of ample angles \cite{R19}, as we hope to discuss elsewhere \cite{R20}.

\subsection{Organization}

Section \ref{flagsec} starts by definining the three types
of blow-ups that arise in this study.
Then, \saldp flags are introduced 
(Definition \ref{badcurvebetadef})
and we state the main result on how such flags characterize strong
asymptotic log positivity (Theorem \ref{badflagcriterionprop}).
To prove this result we derive a characterization of flags containing
a non-boundary component (Lemma \ref{ineqwithdeltaALdPcaseLem}) and those containing
a boundary component (Lemma \ref{bndrycurvesflagslem}). Both of these lemmas are of independent
interest and will serve us repeatedly in the classification of \saldp
flags (Theorem \ref{MainFlagthm}).
In \S\ref{reductionsec} we show any \saldp pair can be
described as a smooth boundary blow-up of another \saldp pair with a smaller
Picard group. The key result here is 
Proposition \ref{noawaynotailProp} 
which is more conceptual approach than that given in \cite[Theorems 2.1, 3.1]{CR}
and which generalizes to the setting of \aldp pairs \cite{CMGR}.
In \S\ref{statesec} we turn to the main application, the 
classification of \saldp pairs (Theorem \ref{mainthm}).
In an Appendix we present an auxiliary classification result (Proposition \ref{baseprop})
under a small Picard rank assumption.

\section{Strongly asymptotically log del Pezzo flags}
\lb{flagsec}

The classification of del Pezzo surfaces is one of the
most classical results in algebraic geometry. It hinges 
on blowing down $-1$-curves and keeping track of intersection
numbers in the process, and in the final step reverse engineering 
to determine the possible location (points) of the blow-downs (blow-ups).

It is therefore natural that when classifying 
strongly asymptotically log del Pezzo we similarly have to study
carefullly the $-1$-curves, this time on the
{\it pair}, i.e., the added layer of 
difficulty is to understand not just the $-1$-curves on the {\it surface }
but also their {\it relative position to the boundary curve.}

Thus, a key step in the proofs of 
Proposition \ref{noawaynotailProp} 
and
Theorem \ref{mainthm} involves distinguishing between three different types of birational operations
on pairs that we introduce in Definition \ref{properdef}.

First, let us establish some basic notation.
 Given a blow-down map 
$$\pi_P:S\ra s$$
of a reduced zero-dimensional locus (i.e., a collection of distinct points)
$$P\subset s$$
and an irreducible curve 
$$\Sigma\subset s,$$
 denote by
$$
\tilde\Sigma\subset S
$$
the $\pi_P$-proper transform of $\Sigma$. For each 
$p\in P$, denote the exceptional curve by
\beq
\lb{excEq}
E_{p}:=\pi^{-1}(p).
\eeq
All log pairs considered will be of the form 
$\big(s,c=\sum_{i=1}^rc_i\big)$ a log pair with each $c_i$ a smooth
irreducible curve in $s$. We will also consider coefficient
vectors $\vec\be=(\be_1,\ldots,\be_r)\in(0,1)^r$ with the convention
\beq
\lb{ciconvEq}
c_0:= c_r, \q \be_0:=\be_r.
\eeq

\bdefin
\lb{properdef}
Let $P,\pi_P$ be as above.

\med
\noi $\bullet\;$
A pair $\big(S,C=\sum_{i=1}^rC_i\big)$ is called a \emph{smooth boundary blow-up} of a pair 
$\big(s,c=\sum_{i=1}^rc_i\big)$ if $P$ is contained in the smooth
locus of $c$ (i.e., $P\cap c_j\cap c_k=\emptyset$ for all $j\not= k\in\{1,\ldots,r\}$)
and $C=\tc$, where $\tc$ is the $\pi_P$-proper transform of $c$.
We say $(s,c)$ is the \emph{smooth boundary  blow-down} of $(S,C)$.

\med
\noi 
$\bullet\;$
A pair $\big(S,C=\sum_{i=1}^rC_i\big)$ is called an \emph{away blow-up} of a pair 
$\big(s,c=\sum_{i=1}^rc_i\big)$ if $p\in s\sm c$ 
and $C=\pi^{-1}(c)=\tc$.
We say $(s,c)$ is the \emph{away blow-down} of $(S,C)$.
We use the same terminology if  $\pi$ is the blow-up of a collection of distinct
points in $s\sm c$.

\med
\noi 
$\bullet\;$
A pair $\big(S,C=\sum_{i=1}^rC_i\big)$ is called a \emph{tail blow-up} of a pair 
$(s,c=\sum_{i=1}^{r-1}c_i)$ if $p\in c_{r-1},\, c_{r-1}.c_i=1$ for precisely one $i\in\{1,\ldots,r-2\}$
and zero otherwise, 
and $C=\pi^{-1}(c)$ (i.e., $C_r=\pi^{-1}(p)$).
We say $(s,c)$ is the \emph{tail blow-down} of $(S,C)$.
We use the same terminology if $\pi$ is the blow-up of a collection of distinct
points  in the smooth locus of $c$,
each on a different tail component.

\edefin

For the purpose of classification of strongly asymptotically
log del Pezzos it turns out that it suffices to consider
just the first type of blow-ups. This is the content 
of
Proposition \ref{noawaynotailProp} 
and
Theorem \ref{mainthm}. 

In order to prove such a result, we establish a characterization of such 
blow-ups that preserve strong asymptotic log positivity
(Theorem \ref{badflagcriterionprop}).
This leads to an additional key new notion we introduce
in this article:

\bdefin
\lb{badcurvebetadef}
Let $P,\pi_P,\Sigma$ be as above.
We say 
$P\subset \Sig 
\subset s$ is
a \emph{\saldp flag (SALdP flag) for $\big(s,c=\sum_{i=1}^rc_i\big)$}
if there exists a sequence
of vectors $\vec\be(j)=(\be_1(j),\ldots,\be_r(j))\in(0,1)^r$ 
tending to the origin such that   
\beq
\lb{badcurvesalfeq}
\bigg(K_{S}+\sum_{i=1}^r(1-\be_i(j))\tilde c_i\bigg).\tilde \Sig\ge0, 
\q \h{for all $j\in\NN$}.
\eeq
\edefin

The point of this definition is the following useful characterization of
strong asymptotic log positivity of a pair in terms of certain SALdP flags.

\bthm
\lb{badflagcriterionprop}
Let $(s,c)$ be a strongly \aldp pair and let $(S,C=\tc)$ be 
the smooth boundary
blow-up at $P=\{p_1,\ldots,p_m\}\subset c$ (Definition \ref{properdef}).
Then $(S,C)$
is strongly \aldp if and only if 
there are no SALdP flags for $(s,c)$ of the form
\beq
\lb{badflagmeq}
\{p_{i_1},\ldots,p_{i_\ell}\}\subset \Sig\subset s,
\eeq
with $\{i_1,\ldots,i_\ell\}\subset\{1,\ldots,m\}$.
\ethm

\bremark
Note that Theorem 
\ref{badflagcriterionprop}
characterizes a positivity property of $(S,C)$
in terms of a `downstairs' pair $(s,c)$ as 
well as curves that live on {\it intermediate} pairs 
$(S',C')$ that are blow-ups of $(s,c)$ but blow-{\it%
downs} of
$(S,C)$! See Remark \ref{DprimeRemark}.
Still, after some thought the content of Theorem \ref{badflagcriterionprop}
might seem intuitive, however, as we will see below there are
a few pitfalls to an `easy' proof. 
\eremark

\subsection{Flags with boundary points and no boundary curve}

In order to prove Theorem \ref{badflagcriterionprop} we will need to develop a basic understanding
of SALdP flags whose 
zero-dimensional locus lies in the boundary.
 The next lemma characterizes such flags under
 the additional assumption that their one-dimensional
locus is {\it not} a boundary component.

\blem
\lb{ineqwithdeltaALdPcaseLem}
Suppose that $\Sig$ is not a component of $c$
and let $P$ be contained in the smooth locus of $c$. Then
$c\supset P=\{p_1,\ldots,p_m\}\subset \Sig\subset s$ is
a SALdP flags for $(s,c)$ if and only if 
\beq
\lb{summqjmiineq}
0\ge 
-(K_s+c).\Sig+ 
\sign(\tc.\tSig),
\eeq
if and only if 
$$
\tc.\tSig=c.\Sig-\sum_{k=1}^m\mult_{p_k}\Sig=0=(K_s+c).\Sig.
$$
\elem

\bpf
First (recall \eqref{excEq}),
\beq
\lb{SigtSigeq}
\tSig
\sim\pi^*\Sig-\sum_{k=1}^m\mult_{p_k}\Sig\, E_{p_k}.
\eeq
Second, as $c\not=0$ and $\{p_k\}_{k=1}^m$ are smooth points of $c$,
\beq
\lb{ttcieq}
\tc
\sim
\pi^*c-\sum_{k=1}^m\mult_{p_k}c\,E_{p_k}
=
\pi^*c-\sum_{k=1}^mE_{p_k}.
\eeq
In addition,
\beq
\lb{KSpiqifirsteq}
K_S
\sim
\pi^*K_s+\sum_{k=1}^mE_{p_k},
\eeq
so
\beq
\lb{KSpiqieq}
K_S+\tc
\sim
\pi^*(K_s+c).
\eeq

The inequalities \er{badcurvesalfeq}, $j\in\NN$, i.e., 
$(K_S+\tc).\tSig\ge \tSig.\sum_{i=1}^r\be_i(j) \tc_i$,
 become
\beq
\lb{summqjmiwithbetasineq}
0\ge -(K_s+c).\Sig+ 
\sum_{\al=1}^r\be_{\al}(j)\Big(c_{\al}.\Sig-\sum_{k=1}^m\mult_{p_k}\Sig\mult_{p_k}c_{\al}\Big), \q\forall j\in\NN.
\eeq
We claim that 
\beq\lb{mpilesigceq}
\sum_{i=1}^m\mult_{p_k}\Sig\mult_{p_k}c_\a\le \Sig.c_\a, \q\forall\al\in\{1,\ldots,r\}.
\eeq
Indeed, $\Sig\not=c_\a$ and thus also $\tSig\not=\tc_\a$, so
since the $p_k$ are smooth points of $c$, by Definition \ref{properdef} and
\eqref{SigtSigeq}--\eqref{ttcieq},
\beq
\lb{tsigtcalphaineq}
0\le \tSig.\tc_\a=\Sig.c_\a-\sum_{k=1}^m\mult_{p_k}\Sig\mult_{p_k}c_\a,
\eeq
as claimed.
Next, by \eqref{SigtSigeq}--\eqref{ttcieq},
$$
\tc.\tSig=
c.\Sig-\sum_{k=1}^m\mult_{p_k}\Sig=
\sum_{\a=1}^r\Big(c_{\al}.\Sig-\sum_{k=1}^m\mult_{p_k}\Sig\mult_{p_k}c_{\al}\Big)
$$
with each summand nonnegative, the sign will equal the sign of the maximal summand.
Finally, both terms in \er{summqjmiwithbetasineq} are nonnegative; the first since $-K_S-C$ is 
nef for every \saldp pair (as a limit of ample classes by Definition \ref{definition:log-Fano}), while the second by \er{tsigtcalphaineq}. Thus both must vanish 
in order for \er{summqjmiwithbetasineq} to hold: the first term is independent of $j$
while the second is either zero or small and positive.
Conversely, if both terms vanish than \er{summqjmiwithbetasineq} holds (with equality)
and so does \er{badcurvesalfeq}, for each $j\in\NN$.
This concludes the proof of the equivalence
of Lemma \ref{ineqwithdeltaALdPcaseLem}.
\epf

\bremark
\lb{crucialseqremark}
The crucial consequence of Lemma \ref{ineqwithdeltaALdPcaseLem} is that
being a SALdP flag with zero-dimensional locus in the boundary and 
one-dimensional locus off the boundary
{\it does not depend at all 
on the range of the $\be_i$'s!} We will use this crucially in the proof below.
This will be generalized and used crucially also in
\cite{R20}.
\eremark

\subsection{Flags with a boundary curve}

 The next lemma characterizes flags whose one-dimensional
locus is a component of the boundary.

\blem
\lb{bndrycurvesflagslem}
Suppose that $\big(s,c=\sum_{\al=1}^rc_\al\big)$ is strongly \aldpno.
Then $\{p_1,\ldots,p_m\}\subset c_i\subset s$ 
is a SALdP flag for $(s,c)$ if and only if:
\beq
\baeq
&\bullet\;
m\ge K_s^2 \h{ if $r=1$ and $c_1\sim-K_s$},
\cr
&\bullet\;
m>c_i^2 \h{ if $r\ge2$ and $c\sim-K_s$},
\cr
&\bullet\;
m>c_i^2 \h{ if $r\ge3$, $c\not\sim-K_s$, and $c_i$ intersects exactly two other $c_j$'s.}
\cr
\eaeq
\eeq
\elem

\bpf
Suppose first that $c_1$ is a smooth elliptic curve so $\mult_{p_k}c_1=1$
for all $k\in\{1,\ldots,m\}$,
$c_1\sim-K_s$, and $r=1$ \cite[Lemma 2.2]{CR}.
Plugging-in $K_s+c\sim0$ in \er{badcurvesalfeq} gives,
$$
\baeq
0
&\ge -(K_s+c).\Sig+
\be_1\Big(c_1^2-\sum_{k=1}^m\mult^2_{p_i}c_1\Big)
=
\be_1(c_1^2-m)
,
\eaeq
$$
i.e., $m\ge c_1^2=K_s^2$.

If $c$ has more than one component 
then each component of $c$ (including $c_i$) must
be a smooth rational curve \cite[\S3]{CR}. 
Let $\eps\in\{0,1,2\}$ be the number of components of $c$
(excluding $c_i$ itself) that $c_i$ intersects, counted with multiplicity. 
By \er{badcurvesalfeq},
\beq
\lb{displayedepsKsSigineq}
\baeq
0
&\ge -(K_s+c).c_i+
\sum_{\al=1}^r\be_{\al}\Big(c_{\al}.c_i-\sum_{k=1}^m\mult_{p_k}c_i\mult_{p_k}c_{\al}\Big)
\cr
&\ge -K_s.c_i-c_i^2-\eps
+
\sum_{\al=1}^r\be_{\al}\Big(c_{\al}.c_i-\sum_{k=1}^m\mult_{p_k}c_i\mult_{p_k}c_{\al}\Big)
\cr
&\ge
2-\eps
+
\sum_{\al=1}^r\be_{\al}\Big(c_{\al}.c_i-\sum_{k=1}^m\mult_{p_k}c_i\mult_{p_k}c_{\al}\Big)
,
\eaeq
\eeq
so it follows that
we must have
$\eps=2$ and, using \cite[Lemma 3.5]{CR}, either 
(i) $r=2$, and $c_i=c_2$ with $c_1.c_2=2$ and $c_1+c_2\sim-K_s$, or
(ii) $r\ge3$, and $c_i$ is a `middle' component of $c$, i.e.,
with $c_i.c_\al=1$ for $\al=i\pm 1\mod r$ 
(recall \er{ciconvEq})
and $c_i.c_\al=0$ for $\al\not\in\{i-1,i,i+1 \mod r\}$  (note that in theory $c$ could have 
several connected components if $c\not\sim-K_s$ \cite[Lemma 3.5]{CR}, but this does not affect the computation).

In case (i), \er{displayedepsKsSigineq} becomes,
$$
\baeq
0
&\ge 
\be_{1}\Big(2-\sum_{k=1}^m\mult_{p_k}c_2\mult_{p_k}c_{1}\Big)
+
\be_{2}\Big(c_2^2-\sum_{k=1}^m(\mult_{p_k}c_2)^2\Big)
\cr
&= 
2\be_{1}
+
\be_{2}(c_2^2-m)
\eaeq
$$
since 
$$
\mult_{p_k}c_2\mult_{p_k}c_{1}=0,
$$
as otherwise $p_k\in c_1\cap c_{2}$
contrary to Definition \ref{properdef}.
Thus,  $c_2$ is part of a SALdP flag
if and only if
$
c_2^2<m,
$
equivalently, $\tc_2^2<0$.

In case (ii), \er{displayedepsKsSigineq} becomes,
$$
\baeq
0
&\ge 
\be_{i-1}\Big(1-\sum_{k=1}^m\mult_{p_k}c_i\mult_{p_k}c_{i-1}\Big)
+
\be_{i}\Big(c_i^2-\sum_{k=1}^m(\mult_{p_k}c_i)^2\Big)
\cr
&\qquad+
\be_{i+1}\Big(1-\sum_{k=1}^m\mult_{p_k}c_i\mult_{p_k}c_{i+1}\Big)
\cr
&= 
\be_{i-1}
+
\be_{i}(c_i^2-m)
+
\be_{i+1}
\cr
\eaeq
$$
since 
$$
\mult_{p_k}c_i\mult_{p_k}c_{i\pm1\!\!\mod r}
=0
$$
(recall \er{ciconvEq}),
as otherwise $p_k\in c_i\cap c_{i\pm1\!\!\mod r}$ 
contrary to Definition \ref{properdef}.
Thus, $c_i$ is part of a SALdP flag  if and only if
$
c_i^2<\sum_{k=1}^m\mult_{p_k}c_i=m,
$
i.e., $\tc_i^2<0$.
\epf

\subsection{Degree under smooth boundary blow-ups}

\blem
\lb{sqtermlemma}
Let $(s,c)$ be \saldp and let $(S,C)$ be obtained from $(s,c)$ via
smooth boundary blow-up. Then $\big(K_S+\sum_{i=1}^r(1-\be_i)C_i\big)^2>0$
for all sufficiently small $\be_i$ if and only if one of the following
three mutually exclusive conditions holds:
\beq
\lb{conditionsforsqrtermlemma}
\baeq
&
\bullet\; m<K_s^2 \h{\ if $c_1\sim-K_s$},
\cr
&
\bullet\; \h{at most $c_i^2$ points are blown on each $c_i$  \h{\ if $c\sim-K_s$ and $r>1$}},
\cr
&
\bullet\; 
c\not\sim-K_s.
\cr
\eaeq
\eeq
\elem

\bremark
For the second condition note that $c_i^2\ge0$ in this case
\cite[Lemmas 3.6, 3.5]{CR}.
\eremark

\bpf
Compute using \er{KSpiqieq},
\begin{equation}
\begin{aligned}
\label{aequivtoEq}
\Big(K_S+\sum_{i=1}^r(1-\be_i)C_i\Big)^2
&=
(K_S+C)^2+\Big(\sum_{i=1}^r\be_iC_i\Big)^2
-
2\sum_{i=1}^r\be_iC_i.(K_S+C)
\cr
&=
(\pi^*K_s+\pi^*c)^2
+\Big(\sum_{i=1}^r\be_iC_i\Big)^2
-2\sum_{i=1}^r\be_ic_i.(K_s+c)
\cr
&=
(K_s+c)^2
+\Big(\sum_{i=1}^r\be_i\tc_i\Big)^2
-2\sum_{i=1}^r\be_ic_i.(K_s+c)
.
\end{aligned}
\end{equation}
There are two possibilities: $c$ is a union of disjoint chains of 
smooth rational curves or else
$c\sim-K_s$ and it is a single cycle \cite[Lemma 3.5]{CR}. 
Observe that $C$ has the same structure
as $c$ by Definition \ref{properdef}.

In the latter case \er{aequivtoEq} reduces to 
$$
\Big(\sum_{i=1}^r\be_i\tc_i\Big)^2>0.
$$
We treat a few sub-cases: if $r=1$ then
\er{aequivtoEq} reduces to $\tc^2>0$. 
If $r=2$, then $c_1.c_2=2$  \cite[Lemma 3.5]{CR} so  
\er{aequivtoEq} reduces to 
$$
\be_1^2\tc_1^2+\be_2^2\tc_2^2+4\be_1\be_2>0,
$$
which holds for all small $\be_i>0$ if and only if $\tc_1^2,\tc_2^2\ge0$.
If $r\ge3$ then $c_i.c_{i-1\!\!\mod r}=c_i.c_{i+1\!\!\mod r}=1$ and otherwise
$c_i.c_j=0$ for $i\not=j$ 
(recall \eqref{ciconvEq}). Thus, 
\er{aequivtoEq} reduces to 
$$
\sum_{i=1}^r\be_i^2\tc_1^2+2\sum_{i=1}^r\be_i\be_{i+1\!\!\!\!\mod r}>0,
$$
which holds for all small $\be_i>0$ if and only if $\tc_i^2\ge0$ for each $i\in\{1,\ldots,r\}$.

In the former case, each $c_i$ is smooth and rational hence $K_s.c_i+c_i^2=-2$. We can assume
without loss of generality that $c$ is connected since the computations of the
linear terms below are done for each connected component. So when
$r=1$ \er{aequivtoEq} reduces to 
$$
(K_s+c)^2
+\be_1^2c_1^2
+4\be_1;
$$
when $r=2$ we obtain 
$$
(K_s+c)^2
+\big(\sum_{i=1}^2\be_ic_i\big)^2
+2\be_1+2\be_2;
$$
when $r>2$ we obtain
$$
(K_s+c)^2
+\big(\sum_{i=1}^r\be_ic_i\big)^2
+2\be_1+2\be_r;
$$
all three of these expressions are positive for all small $\be_i>0$
as the linear terms have positive coefficients (and the constant term
is nonnegative since $-K_s-c$ is nef (recall the end of the proof of Lemma \ref{ineqwithdeltaALdPcaseLem}). 
Here we used that
 $c_1.c_2=c_i.c_{i+1}=\ldots
=c_{r-1}.c_r=1$ for all $i=1,\ldots,r-1$ with all other $c_i.c_j=0$
for $i\not=j$. 
\epf


\subsection{Proof of Theorem \ref{badflagcriterionprop}}

Let us first assume that $(S,\tc)$ is strongly \aldpno.
If there exists a SALdP flag of the form \er{badflagmeq} then \er{badcurvesalfeq}
and the Nakai--Moishezon criterion \cite[Theorem 1.2.23]{Laz} imply that
that $(S,\tc)$ is not \saldp by Definition \ref{definition:log-Fano},
a contradiction.

Suppose now that no SALdP flags of the form \er{badflagmeq} exist. 
If $(S,\tc)$ is not \saldp there exists a sequence $\{\be(j)\}_{j\in\NN}$
tending to the origin in $\RR^r$ such that the class
$\big(K_{S}+\sum_{i=1}^r(1-\be_i(j))\tilde c_i\big)$ is not negative
for each $j$.
By Nakai--Moishezon this means that either 
\beq\lb{NMnotholdEq}
\bigg(K_{S}+\sum_{i=1}^r(1-\be_i(j))\tilde c_i\bigg)^2
\le 0, \q \h{for all $j\in\NN$},
\eeq
and/or there exists {\it a sequence} of irreducible curves $0\not=D_j\subset S$
such that 
\beq
\lb{NM2ndnotholdEq}
\bigg(K_{S}+\sum_{i=1}^r(1-\be_i(j))\tilde c_i\bigg).D_j\ge 0,\q
\h{for all $j\in\NN$}.
\eeq
The first possibility \er{NMnotholdEq} cannot hold by 
Lemma \ref{sqtermlemma} 
{\it if we can verify the conditions \er{conditionsforsqrtermlemma}}.
This is the first subtle point in the proof of Theorem \ref{badflagcriterionprop}. Fortunately,
our assumption that no SALdP flags of the form \er{badflagmeq} exist precisely
verifies these conditions: Lemma \ref{bndrycurvesflagslem} shows that if any of
these conditions is not verified then one of the boundary components of 
$c$ would be the one-dimensional strata of a SALdP flag of the form \er{badflagmeq} (this is only needed in 
the case $c\sim-K_s$)!

Thus, only the second possibility, i.e., \er{NM2ndnotholdEq}, can hold.
A subtle, but important, point is that
we can assume \er{NM2ndnotholdEq} holds with a {\it fixed} $D=D_j$:

\bclaim
\lb{fixedDClaim}
There exists a fixed irreducible curve $D\subset S$ such that
$$
\bigg(K_{S}+\sum_{i=1}^r(1-\be_i(j))\tilde c_i\bigg).D\ge 0,\q
\h{for all $j\in J\subset\NN$ with $J$ an infinite set}.
$$
\eclaim
\bpf The proof is based on Remark \ref{crucialseqremark}, 
but requires a bit more. 

Let $\calD=\{D_j\,:\, j\in\NN\}$ be the collection of divisors
satisfying \er{NM2ndnotholdEq}. If infinitely-many of the $D_j$
are boundary components, then there exists an unbounded subsequence $J\subset\NN$ 
with $D_j=\tilde c_\al=C_\al$ for some fixed $\al$
and for all $j\in J$,
and we are done. Otherwise, let $J:=\{j\in\NN\,:\, D_j\not\subset C\}$.
By Remark \ref{crucialseqremark} we may take $D=D_j$ for an arbitrary
(fixed) $j\in J$!
\epf

We return to the proof of Theorem \ref{badflagcriterionprop}.
Let $D$ be the curve furnished by Claim \ref{fixedDClaim}.
Write 
$$
D=\pi^*\pi(D)-\sum_{k=1}^m\mult_{p_k}\pi(D)E_{p_k},
$$
so using \er{ttcieq}, \er{KSpiqieq}, and Claim \ref{fixedDClaim} there exists 
an unbounded subsequence $J\subset\NN$
 such that for all $j\in J$,
\beq
\lb{DineqonSeq}
\baeq
0
&\le 
\Big(K_{S}+\sum_{i=1}^r(1-\be_i(j))\tilde c_i\Big).D
\cr
&=
\Big(K_{s}+\sum_{i=1}^r(1-\be_i(j)) c_i\Big).\pi(D)
+\sum_{i=1}^r\be_i(j)\sum_{k=1}^m\mult_{p_k}\pi(D)\mult_{p_k}c_i.
\cr
\eaeq
\eeq
Since
$(s,c)$ is \saldp the first term in the last line is negative for all sufficiently
large $j\in J$.
By assumption, each $p_k$ is a smooth point of $c$, i.e.,
$\mult_{p_k}c_i=1$ for some (actually, exactly one) $i\in\{1,\ldots,r\}$ it follows
that $\mult_{p_k}\pi(D)>0$ for at least one $k\in\{1,\ldots,m\}$.
Thus, by Definition \ref{badcurvebetadef},
$$
\emptyset\not=\{p_k\,:\, \mult_{p_k}\pi(D)>0\}\subset\pi(D)\subset s
$$
is a SALdP flag for $(s,c)$ of the form \er{badflagmeq}, contradicting
our assumption. Thus, $(S,\tc)$ is \saldpno, concluding the proof
of Theorem  \ref{badflagcriterionprop}.

\bremark
\lb{DprimeRemark}
Let $\pi':S'\ra s$ be the blow-up of $s$  at the subset of points
$
\{p_k\,:\, \mult_{p_k}\pi(D)>0\}\subset\pi(D)\subset s.
$
Then $D':=\widetilde{\pi(D)}^{\pi'}$ lives on the 
`intermediate' surface $S'$ `between' $s$ and $S$
(unless of course $\pi(D)$ passes through all $p_1,\ldots,p_m$).
A moment's thought shows that $D'$ satisfies the same
exact inequality \er{DineqonSeq} $D$ satisfied on $S$, but
on $S'$ instead. Thus, we do get a SALdP flag for $(s,c)$ 
as claimed.
\eremark


\section{Reduction to smooth boundary blow-ups}
\lb{reductionsec}

\subsection{Trichotomy of $-1$-curves}

By a {\it $-1$-curve} we will always mean an irreducible smooth rational curve
of self-intersection $-1$ (and genus $g$ zero). 
Let $e$ be a $-1$-curve on an asymptotically log del Pezzo 
pair $(s,c=\sum_{i=1}^rc_i)$. By adjunction \cite{GH}, 
\beq\lb{adjEq}
K_s.e=2g(e)-2-e^2=-1,
\eeq 
so
$$
0> \big(K_s+\sum_{i=1}^r(1-\be_i)c_i\big).e
=-1+c.e+O(\vec\be),
$$
i.e., $c.e<2$. 
If $e$ is not contained in $c$, i.e., $e\not=c_i$ for any $i$, 
then also $c.e\ge0$.
Thus, on any \aldp pair $(s,c)$ there are 
three disjoint families of $-1$-curves:

\med
\noi 
$$
\baeq
&\bullet\; 
\calEprop(s,c)
:=
\{E\,:\,E\not\subset C \h{ is a $-1$-curve with $E.c=1$}\}
,%
\cr
&\bullet\; \calEaway(s,c)
:=
\{E\,:\,E\not\subset c \h{ is a $-1$-curve with $E.c=0$}\}
,
\cr
&\bullet\; \calEC(s,c) 
:=
\{E\,:\,E\subset c \h{ is a $-1$-curve}\},
\eaeq
$$
and the disjoint union 
\beq
\lb{calEeq}
\calE(s):=\calEprop(s,c)
\cup\calEaway(s,c)\cup\calE_{C}(s,c)
\eeq
consists of {\it all} $-1$-curves on $s$ \cite[Lemma 3.3]{CR}.
What is more, by Remark \ref{alltailsremark} below, this trichotomy precisely corresponds to the three birational operations in Definition \ref{properdef}.

\bremark
\lb{alltailsremark}
In fact, $\calE_{C}(s,c)$ consists exclusively of tail $-1$-components
of the boundary \cite[Lemma 3.6]{CR}.
\eremark

\subsection{Smooth boundary blow-ups suffice}

\begin{figure}
		\centering

\tikzstyle{block} = [rectangle, draw, 
    text width=7.5em, text centered, rounded corners, minimum height=3.25em]
\tikzstyle{line} = [draw, -latex']

\begin{tikzpicture}[node distance = 4.05cm, auto] 
    \node [block] (init) {$C\sim-K_S$?};
    \node [block, left of=init] (expert) {contract $\calEprop(S,C)$ to obtain $(S(1),C(1))$};
    \node [block, right of=init] (system) {$\rk\Pic(S(1))\le2$};
    \node [block, below of=init] (identify) {$\calEC(S(1),C(1))$ empty?}; 
    \node [block, left of=identify] (evaluate) {contract $\calEaway(S,C)$ to obtain $(S(2),C(2))$};
    \node [block, below of=evaluate] (decide) {$r=1$?};
    \node [block, right of=decide] (update) {contract a subset of $\calEC(S(2),C(2))$ to obtain $(S(3),C(3))$ with $\rk\Pic(S(3))=2$};
    \node [block, below of=decide] (stop) {$\calEprop(S(1),C(1))$ nonempty, contradiction}; 
    \node [block, right of=update] (new) {$\calEprop(S(2),C(2))$ nonempty, contradiction}; 
    \path [line] (init) -- node {no} (identify);
    \path [line] (identify) -- node {no} (evaluate);
    \path [line] (evaluate) -- (decide);
    \path [line] (decide) -- node {no} (update);
    \path [line] (identify) -- node {yes} (system);
    \path [line] (decide) -- node {yes}(stop);
    \path [line] (expert) --  (init);
    \path [line] (init) -- node {yes} (system);
    \path [line] (update) -- (new);
\end{tikzpicture}

		\caption{A flowchart for the proof of Proposition \ref{noawaynotailProp}.}
		\label{Figure0}
	\end{figure}
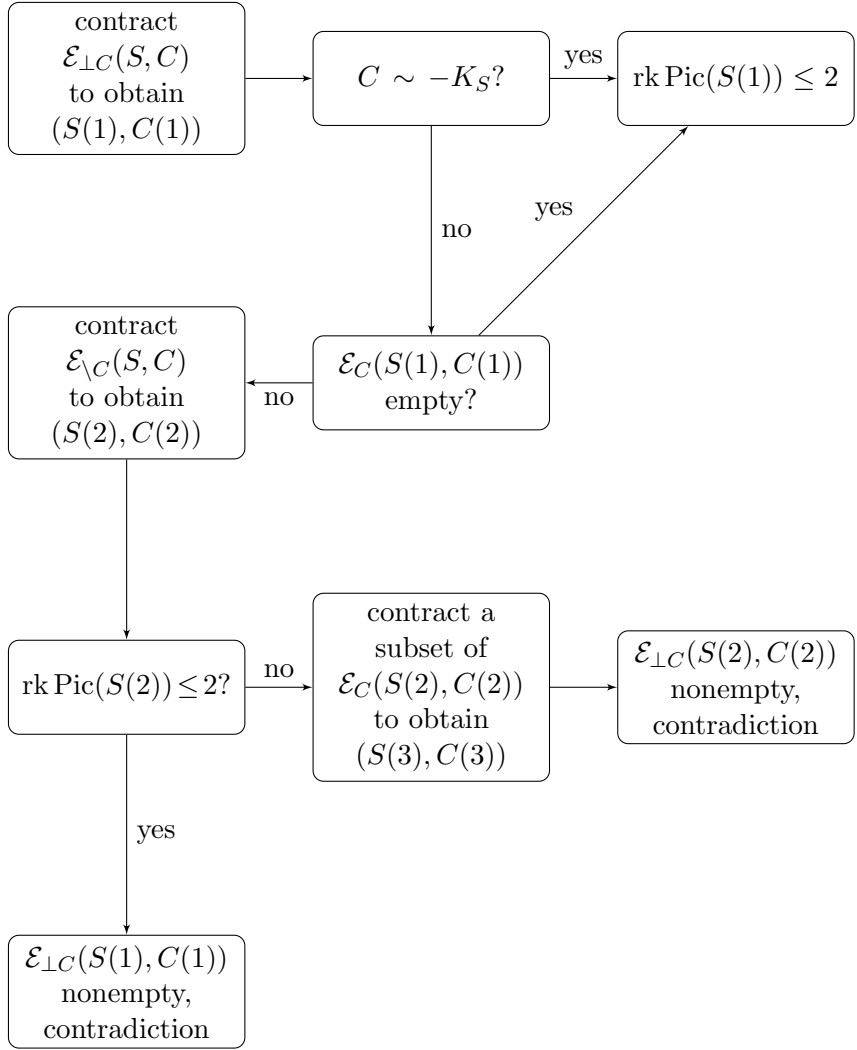

The next result is key to the proof of Theorem \ref{mainthm}. 
The proof is modelled on the one in \cite{CR} but is more
conceptual.
It shows
that the first type of blow-ups 
in Definition \ref{properdef}
suffice as far as the classification
of \saldps is concerned.

\bprop
\lb{noawaynotailProp} 
Let $(S,C)$ be a \saldp pair. Let
$(s,c)$ be the smooth boundary  blow-down of $(S,C)$
given by the contraction of $\calEprop(S,C)$.
Then $(s,c)$ is a \saldp pair with $\rk\Pic(s)\le 2$.
\eprop

\bpf
A flowchart for the proof is provided in Figure \ref{Figure0}.
Note that $S$ and hence $s$ are rational surfaces \cite[p. 1253]{CR}.
We will use repeatedly the following consequences of the classification
of rational surfaces \cite[p. 520]{GH}:
\beq\label{MinusOneCurveExistence}
\h{every rational surface with $\h{rk}(\h{Pic})>2$ contains a $-1$-curve},
\eeq
and that
\beq
\label{RationalClassifPicAtMostTwo}
\h{a rational surface with $\h{rk}(\h{Pic})\le 2$
is either $\PP^2$ or $\FF_n, \,n\ge0$}.
\eeq

\noi 
{\it Step 1: the case $C\sim-K_S$. } Suppose first that $C\sim-K_S$. Then $\calE_{C}(S,C)=\emptyset$
\cite[Lemma 3.6]{CR} and $\calEaway(S,C)=\emptyset$
since $-K_S.E=1$ by \eqref{adjEq}, and this equals $C.E$. 
So by \eqref{calEeq}, $\calE(S)=\calEprop(S,C)$
and contracting all of these curves is a smooth boundary  blow-down that
yields a pair $(s,c)$ with no $-1$-curves. 
 By
 \eqref{MinusOneCurveExistence} then
$\rk(\Pic(s))\le2$ and, finally, the pair $(s,c)$ is \saldp
\cite[Lemma 3.4]{CR}, as desired.

\medskip\noi 
{\it Step 2: contracting $\calEprop(S,C)$.  }
From now on, suppose $C\not\sim-K_S$.
Let 
$$
\pi_1:S\ra S(1)
$$
be the contraction of $\calEprop(S,C)$ 
and let 
$$
(S(1),C(1))
$$ 
be the smooth boundary  blow-down of $(S,C)$.
The pair $(S(1),C(1))$ is \saldp as $(S,C)$ is by assumption
\cite[Lemma 3.4]{CR}.
\bclaim
\lb{firstClaimSone}
Assume $C\not\sim-K_S$ and let 
$
(S(1),C(1))
$
be the smooth boundary  blow-down of $(S,C)$.
Then,
\beq
\lb{noperpSone}
\calEprop(S(1),C(1))=\emptyset.
\eeq
\eclaim
\bpf
Denote $\{E_1,\ldots,E_m\}=\calEprop(S,C)$
and $q_i:=\pi_1(E_i)\in C(1)$ (the last inclusion
is since otherwise $E_i$ would not intersect
$C$ but by assumption $E_i.C=1$).
Let $e\in \calEprop(S(1),C(1))$. 
Then $\tilde e\subset S$ satisfies
$\tilde e^2=(\pi_1^*e-\sum_{i=1}^m\mult_{q_i}e\,E_i)^2\le e^2=-1$.
Equality is not possible since then $\tilde e\in \calEprop(S,C)$
and would be contracted by $\pi_1$ to a point and not to a curve
$e$. So $\tilde e^2\le -2$. On the other hand, $e\not\subset C(1)$
so  $\tilde e\not\subset C$ and (compare to the proof of
the case $r=1$ \cite[Lemma 2.5]{CR})
\beq
\lb{Lemma2.5SNCEq}
\begin{aligned}
\tilde e^2 &=2h^1(\calO_{\tilde e})-2-K_S.\tilde e\cr
&\ge -2-\big(K_S+\sum_{i=1}^r(1-\be_i)C_i\big).\tilde e
+
\sum_{i=1}^r(1-\be_i)C_i.\tilde e\cr
&\ge-2-\big(K_S+\sum_{i=1}^r(1-\be_i)C_i\big).\tilde e
>-2.
\end{aligned}
\eeq
Therefore, $\tilde e$ cannot exist on $S$, so neither can $e$ on 
$S(1)$.
\epf

\medskip\noi 
{\it Step 3: away blow-down of  $(S(1),C(1))$.  }
Let 
\beq
\lb{pi2mapeq}
\pi_2:S(1)\ra S(2)
\eeq
be the
blow-down map of $\calEaway(S(1),C(1))$ and denote by
$$
(S(2),C(2))
$$
the away blow-down of $(S(1),C(1))$.
The pair $(S(2),C(2))$ is \saldp as 
$(S(1),C(1))$ is (the latter was noted right before Claim \ref{firstClaimSone})
\cite[Lemma 3.4]{CR}.
Observe that
\beq
\lb{awayStwoEq}
\calEaway(S(2),C(2))=\emptyset.
\eeq
Indeed, by the same reasoning as in the proof of Claim \ref{firstClaimSone}
(see \eqref{Lemma2.5SNCEq}),
if there were a $-1$-curve disjoint from $C(2)$ then its
$\pi_2$-proper transform would be a curve disjoint from $C(1)$
and of self-intersection
at most, and hence exactly, $-1$, but then it would be contracted
by $\pi_2$.  
Next, observe that
by \er{noperpSone} and Claim \ref{nonewcurvesClaim} below also
\beq\lb{perpStwoEq}
\calEprop(S(2),C(2))=\emptyset.
\eeq
As $\pi_2$ is an isomorphism on $C(2)$, i.e., $C(1)$ and $C(2)$
are isomorphic and $C(1)_i\cong C(2)_i$
for each $i\in\{1,\ldots,r\}$, the self-intersection numbers
of the boundary components are unchanged under $\pi_2$: 
$C(1)_i^2=C(2)_i^2$. Thus,
\beq\lb{CStwoEq}
\calE_{C}(S(2),C(2))\cong
\calE_{C}(S(1),C(1)).
\eeq

\medskip\noi 
{\it Step 4: the case $C(1)$ contains no $-1$-curves.  }

\bclaim
\lb{awayCClaim}
Assume $C\not\sim-K_S$ and let 
$
(S(1),C(1))
$
be the smooth boundary  blow-down of $(S,C)$.
If
$\calE_{C}(S(1),C(1))=\emptyset$
then either
$\calEaway(S(1),C(1))=\emptyset$
or\hfill\break
$(S(1),C(1))\in \{$\hbox{\rm (I.3A), (I.3B), (II.3)}$\}$.
In particular, $\rk\Pic(S(1))\le 2$.
\eclaim

\bpf
Suppose that $\calE_{C}(S(1),C(1))=\emptyset$.
By \eqref{CStwoEq} also $\calE_{C}(S(2),C(2))=\emptyset$.
So by \eqref{awayStwoEq}--\eqref{perpStwoEq} and \er{calEeq},
$\calE(S(2))=\emptyset$ and
by \eqref{MinusOneCurveExistence}
then $\rk\Pic(S(2))\le 2$. 
Additionally, 
$C(2)\not\sim-K_{S(2)}$:
by 
\eqref{KSpiqieq},
$$
-K_{S(1)}-C(1)
\sim
-\pi^*(K_{S(2)}+C(2))-\sum_{i=1}^kE_{i},
$$
and the left hand side is a limit of ample
classes, while if 
$C(2)\sim-K_{S(2)}$, the right hand side
would have negative square.
In sum:

\smallskip
\noi
$\bullet\;$ $\calE(S(2))=\emptyset$,

\noi
$\bullet\;$ $\rk\Pic(S(2))\le 2$,

\noi
$\bullet\;$ $C(2)\not \sim-K_{S(2)}$. 

\noi
$\bullet\;$ $(S(2),C(2))$ is \saldpno.

\smallskip
\noi
By Proposition 
\ref{baseprop}, $(S(2),C(2))$
must be one of the following 11 pairs:
\beq
\begin{aligned}
\lb{elevenpairsEq}
&\hbox{(I.1B), (I.1C), (II.1B),}\cr
&\hbox{(I.2.n),  (I.4B), (I.4C), (II.2A.n), (II.2B.n), 
(II.2C.n), (II.4B), or (III.3.n).}
\end{aligned}
\eeq
There are two cases to consider. 

First, each of the last 8 surfaces is ruled and contains a fiber through every point.
If exists $T\in \calEaway(S(1),C(1))$, 
let $F$ be the fiber (i.e., $F^2=0$) through the point $\pi_2(T)\not\in C(2)$. In each
of the above cases $F.C(2)\ge1$. More precisely, one can arrange (i.e., in the case
(I.4B) we may choose $F$ to belong to either fibration resulting in different
intersection numbers):
$$
F.C(2)=
\begin{cases}
1,&
\hbox{$(S(2),C(2))\in\{$\rm (I.2.n),  (I.4B), (I.4C), 
(II.2C.n)$\}$},
\cr
2&
\hbox{$(S(2),C(2))\in\{$\rm (I.4B), (II.2A.n), (II.2B.n), 
(II.4B), (III.3.n)$\}$}.
\cr
\end{cases}
$$

The $\pi_2$-proper transforms of $F$ and of $C(2)$ are curves on $S(1)$.
Note that
$\tilde C(2)=\pi^*C(2)=C(1)$. Since $(S(1),C(1))$ is \saldpno, 
the curve $\tilde F$ on $S(1)$ must satisfy $\tilde F^2\ge -1$
by the argument of \eqref{Lemma2.5SNCEq}
(cf. \cite[Lemma 2.5]{CR} for the case $r=1$). So while by choice $\pi_2(T)\in F$, no other
blow-up point of $\pi_2$ can contained in $F$. Thus, 
$\tilde F=\pi_2^*F-T,$ and $\tilde F^2=-1$. Moreover, since
by assumption $T.C(1)=0$,
$$
F.C(2)=(\pi_2^*F-T).(\pi_2^*C(2))=\tilde F. \tilde C(2)=\tilde F.C(1).
$$
Thus, the possibility $F.C(2)=1$ is impossible as then $\tilde F\in 
\calEprop(S(1),C(1))$ contradicting
\eqref{noperpSone}.
The possibility $F.C(2)=2$ is also impossible as then $\tilde F$
is a $-1$-curve intersecting $C(1)$ at 2 points, contradicting 
the fact that $(S(1),C(1))$ is \saldp (recall, once again, \eqref{Lemma2.5SNCEq}).
Altogether, we have shown that if $(S(2),C(2))$ 
is one of the last 8 pairs in \eqref{elevenpairsEq} then $T$
cannot exist, i.e., 
$\calEaway(S(1),C(1))=\emptyset$. 

Suppose now that $(S(2),C(2))$ 
is one of the first 3 cases in \eqref{elevenpairsEq}, namely, 
$(S(2),C(2))\in\{$(I.1B), (I.1C), (II.1B)$\}$, so $S(2)=\PP^2$.
If $\calEaway(S(1),C(1))\not=\emptyset$, denote
$$\pi_2\big(\calEaway(S(1),C(1))\big)=\{q_1,\ldots,q_k\},$$ 
a nonempty collection of points
in $S(2)=\PP^2$ away from $C(2)$. Blowing up any one of these points,
say $q_1$,
we obtain a map $\hat\pi_2:\hat S(2)\ra S(2)$.
Note that 
 $\Big(\hat S(2),\hat C(2)=\hat\pi_2^{-1}\big(C(2)\big)\Big)$ is \saldp as it
is an away blow-down of 
$(S(1),C(1))$ at $\pi_2^{-1}(q_2),\ldots,$\hfill\break $\pi_2^{-1}(q_k)$
\cite[Lemma 3.4]{CR}.
If $k=1$ then $S(1)= \hat S(2)=\FF_1$ and
$(S(1),C(1))\in\{$(I.3A), (I.3B), (II.3)$\}$.
The cases $k\ge2$ are not possible by the same analysis of the previous paragraph
(i.e., there is a fiber on $\hat S(2)$ whose proper transform 
on $S(1)$ would either be an element of 
$\calEprop(S(1),C(1))$ (contradicting
\eqref{noperpSone}) or else a $-1$-curve intersecting the boundary more than once
(contradicting the argument of \eqref{Lemma2.5SNCEq})).
Thus, we have shown
that $(S(1),C(1))\in \{$(I.3A), (I.3B), (II.3)$\}$.

In sum, 
either
$\calEaway(S(1),C(1))=\emptyset$
or
$(S(1),C(1))\in \{$\hbox{\rm (I.3A), (I.3B), (II.3)}$\}$, and in
both cases $\rk\Pic(S(1))\le 2$ (recall  \eqref{MinusOneCurveExistence}),
as stated in Claim \ref{awayCClaim}.
\epf

\bclaim
\lb{nonewcurvesClaim}
Let $(M,A)$ be a \saldp pair with $\calE_{\perp C}(M,A)=\emptyset$
and let $(M',A')$ be the pair obtained by a composition of any number of away
and/or tail blow-downs (Definition \ref{properdef}). Then 
$\calE_{\perp C}(M',A')=\emptyset$.
\eclaim
\bpf
By the same reasoning as in the proof of Claim \ref{firstClaimSone},
any blow-down map can only increase the self-intersection
of curves. However,  there are no $-n$-curves in $M$
intersecting $A$ transversally: $n=1$ by assumption and
$n\ge2$ by \cite[Lemma 2.5]{CR}.
Thus, no new $-1$-curves intersecting $A'$ can appear downstairs.
\epf

\medskip
\noi 
{\it Step 5: the case $C(1)$ contains a $-1$-curve and $r=1$.   }
Suppose that $C(1)$ contains a $-1$-curve and $r=1$. Then
 $\calE_{C}(S(1),C(1))=\{C(1)\}$. By \eqref{CStwoEq},
$\calE_{C}(S(2),C(2))=\{C(2)\}$.
Together with \eqref{awayStwoEq}--\eqref{perpStwoEq} and \er{calEeq},
$\calE(S(2))=\{C(2)\}$.
By \eqref{MinusOneCurveExistence} then $\rk\Pic(S(2))\le 2$
 and in fact, by Proposition
\ref{baseprop}, $(S(2),C(2))$
must be (I.2.1) (recall from the proof of Claim \ref{awayCClaim}
that $(S(2),C(2))$ is \saldpno). In particular, 
 $S(2)=\FF_1$ using the notation of the proof of Claim \ref{awayCClaim}
 the $\pi_2$-proper transform of a fiber through the point $q_1$
 is an element of 
$\calE_{\perp C}(S(1),C(1))$, contradicting \eqref{noperpSone}.

\medskip
\noi 
{\it Step 6: the case $C(1)$ contains a $-1$-curve and $r>1$.   }
Suppose that $C(1)$ contains a $-1$-curve and $r>1$. If
 $\calE_{C}(S(1),C(1))$ is a singleton we are reduced to Step 5
above. So assume that $C(1)$, and hence also $C(2)$, contains 
at least two $-1$-curves. Together with \eqref{RationalClassifPicAtMostTwo} 
this implies that $\rk\Pic(S(2))> 2$.

Note that by Remark \ref{alltailsremark}
together with \eqref{awayStwoEq}--\eqref{perpStwoEq} and \er{calEeq},
all $-1$-curves on $S(2)$
are are mutually disjoint tail components of $C(2)$.
Consider the tail blow-down
$$
\pi_3:S(2)\ra S(3)
$$ 
of a nonempty subset of $\calE_{C}(S(2),C(2))=\calE(S(2))$
so that $\rk\Pic(S(3))=2$.
Denote by
$$
(S(3),C(3))
$$
the resulting tail blow-down of $(S(2),C(2))$.
Note $(S(3),C(3))$ is \saldp \cite[Lemma 3.12]{CR}.
By Proposition 
\ref{baseprop} $(S(3),C(3))$ is one of the following
(recall that $-K_{S(3)}\not\sim C(3)$ as we showed earlier
that
$-K_{S(2)}\not\sim C(2)$):
(I.2.n), (I.3A), (I.3B), 
(I.4B), (I.4C), (II.2A.n), (II.2B.n), 
(II.2C.n), (II.3), 
or (III.3.n).
Each of these surfaces are ruled and contain a fiber through every point.
Let $F$ be the fiber through the point $\pi_3(T)\in C(3)$ with $T$
in the exceptional locus of $\pi_3$ (and hence $T\in
\calE_{C}(S(2),C(2))$). Observe that $\pi_3(T)$ is a smooth point of $C(3)$
so it cannot be the intersection point of the two components of (II.2C.n).
The $\pi_3$-proper transform of $F$ is then a $-1$-curve in $S(2)$ not contained
in $C(2)$ but intersecting $C(2)$ at $T\in C(2)$ transversally,
i.e., $\tilde F\in \calEprop(S(2),C(2))$. This contradicts
\eqref{perpStwoEq}. In conclusion then $T$ cannot exist, i.e.,
$\calE_{C}(S(2),C(2))=\calE(S(2))$ must both be empty.
This means that also that $\calE_{C}(S(1),C(1))=\emptyset$, i.e.,
the hypothetical case of Step 6 is impossible.

Altogether, combining the 6 steps above, the proof of 
of 
\PR{noawaynotailProp} is complete.
\epf

\section{Classification of SALdP flags with zero-dimensional
locus in the boundary}

Motivated by Proposition \ref{noawaynotailProp}, we restrict our
attention in this article to smooth boundary blow-ups.
The next result builds on the tools of \S\ref{flagsec} to completely classify
SALdP flags arising from blow-ups of points on the boundary.

\bthm
\lb{MainFlagthm}
Let $\big(s,c=\sum_{i=1}^{r}c_i\big)$ be \saldp pair that is not the
smooth boundary 
blow-up  (recall Definition \ref{properdef}) of any other \saldp pair. Suppose also that
$\rk(\Pic(s))\le 2$ so that $(s,c)$ is one of the pairs listed 
in \PR{baseprop}.
	Then $P_\pi\subset\Sig\subset s$ is a SALdP flag with $P_\pi\subset c$ if and only if it
	is one of the following (the numbering of the pairs matches
	the notation of \PR{baseprop}):

\begin{itemize}
\item [$\mathrm{(I.1A})$] 
$c\cap \Sig =\{p_1,\ldots,p_{m}\}\subset \Sig\not=c$ with $m>0$, 
so that $\sum_{i=1}^m\mult_{p_k}\Sig=\Sig.c$, i.e., $\tSig.\tc=0$,

\noi
$\{p_1,\ldots,p_9,\ldots,p_{m}\}\subset \Sig=c$ with $ \, m\ge9$,

\item [$\mathrm{(I.3A})$] 
$
c\cap \Sig=\{p_1,p_2\}
\subset \Sig=\h{fiber}, 
$

\item [$\mathrm{(I.4A})$] 
$c\cap \Sig =\{p_1,\ldots,p_{m}\}\subset \Sig\not=c$ with $m>0$, 
so that $\sum_{i=1}^m\mult_{p_k}\Sig=\Sig.c$, i.e., $\tSig.\tc=0$, 

\noi
$\{p_1,\ldots,p_8,\ldots,p_{m}\}\subset \Sig=c$ with $ \, m\ge8$,

\item [$\mathrm{(I.4B})$] 
$
c\cap \Sig=\{p_1,p_2\}\subset \Sig=\h{$(0,1)$-curve},
$

\item [$\mathrm{(II.1A})$] 
$c\cap \Sig =\{p_1,\ldots,p_{m}\}\subset \Sig\not=c$ with $m>0$, 
so that $\sum_{i=1}^m\mult_{p_k}\Sig=\Sig.c$, i.e., $\tSig.\tc=0$, 

\noi
$\{p_1,p_2,p_3,p_4,p_5,\ldots,p_{m}\}\subset c_1$ with $ \, m\ge5$,

\noindent
$\{p_1,p_2,\ldots,p_{m}\}\subset c_2$ with $ \, m\ge2$,

\item [$\mathrm{(II.2A.n})$]
$
c\cap \Sig=\{p_1,p_2\}
\subset \Sig=\h{fiber}, 
$

\item [$\mathrm{(II.2B.n})$]
$
c\cap \Sig=\{p_1,p_2\}
\subset \Sig=\h{fiber},
$

\item [$\mathrm{(II.3})$]
$
c\cap \Sig=\{p_1,p_2\}
\subset \Sig=\h{fiber},
$

\item [$\mathrm{(II.4A})$]
\noindent
$c\cap \Sig =\{p_1,\ldots,p_{m}\}\subset \Sig\not=c_i$ with $m>0$, 
so that $\sum_{i=1}^m\mult_{p_k}\Sig=\Sig.c$, i.e., $\tSig.\tc=0$,

\noindent
$\{p_1,p_2,p_3,\ldots,p_{m}\}\subset 
c_i, i\in\{1,2\}$ with $\, m\ge3$, 

\item [$\mathrm{(II.4B})$]
$
c\cap \Sig=\{p_1,p_2\}
\subset \Sig=\h{$(0,1)$-fiber},
$

\noindent
$\{p_1,p_2,p_3,p_4,p_5,\ldots,p_{m}\}\subset 
c_1$, with $\, m\ge5$, 

\noindent
$\{p_1,\ldots,p_{m}\}\subset 
c_2$, with $\, m\ge1$,

\item [$\mathrm{(III.1})$]
$c\cap \Sig =\{p_1,\ldots,p_{m}\}\subset \Sig\not=c_i$ with $m>0$, 
so that $\sum_{i=1}^m\mult_{p_k}\Sig=\Sig.c$, i.e., $\tSig.\tc=0$,

\noindent
$\{p_1,p_2,\ldots,p_{m}\}\subset c_i, i\in\{1,2,3\}$ with $ \, m\ge2$,

\item [$\mathrm{(III.2})$] 
$c\cap \Sig \supset\{p_1,\ldots,p_{m}\}\subset \Sig\not=c_i$ with $m>0$, 
so that $\sum_{i=1}^m\mult_{p_k}\Sig=\Sig.c$, i.e., $\tSig.\tc=0$,

\noindent
$\{p_1,p_2,p_3,\ldots,p_{m}\}\subset c_1$, with $ \, m\ge3$, 

\noindent
$\{p_1,\ldots,p_{m}\}\subset c_i, i\in\{2,3\}$ with $\, m\ge1$, 

\item [$\mathrm{(III.3.n})$] 
$
c\cap \Sig=\{p_1,p_2\}
\subset \Sig=\h{fiber}\not=c_2,
$

\noindent
$\{p_1,\ldots,p_{m}\}\subset c_2$ with $\, m\ge1$,

\item [$\mathrm{(IV})$]	
$c\cap \Sig \supset\{p_1,\ldots,p_{m}\}\subset \Sig\not=c_i$ with $m>0$, 
so that $\sum_{i=1}^m\mult_{p_k}\Sig=\Sig.c$, i.e., $\tSig.\tc=0$,

\noindent
$\{p_1,\ldots,p_{m}\}\subset c_i, i\in\{1,2,3,4\},$ $\, m\ge1$, 
\end{itemize}

\ethm

\bpf
\PR{baseprop} gives the list of \saldp pairs
 $(s,c)$ with $\rk(\Pic(s))\le 2$ that are not the smooth boundary
blow-up (Definition \ref{properdef}) of any other \saldp pair. 
Applying to this list case-by-case Lemma \ref{ineqwithdeltaALdPcaseLem} readily gives all
SALdP flags $P_\pi=c\cap\Sig\subset\Sig\not\subset c$ while
Lemma \ref{bndrycurvesflagslem}
gives all SALdP flags $P_\pi\subset\Sig\subset c$.
\epf

\section{Classification of \saldp pairs}
\lb{statesec}

The following theorem describes all \saldp pairs
in terms of smooth boundary blow-ups (recall Definition \ref{properdef}).
It is due to Cheltsov and the author \cite[Theorems 2.1, 3.1]{CR}.

	\begin{theorem}
		\label{mainthm} Let $S$ be a smooth surface, let
		$C_1,\ldots,C_r$ be distinct irreducible smooth curves on $S$ such that
		$\sum_{i=1}^{r}C_i$  is a divisor with simple normal crossings.
		Then $\big(S,\sum_{i=1}^{r}C_i\big)$ is strongly
		asymptotically log del Pezzo if and only if it
		is one of the following pairs:
		\begin{itemize}
			\item [$\mathrm{(I.1A})$] $S=\mathbb{P}^2$, $C_1$ is a cubic,%
			\item [$\mathrm{(I.1B})$] $S=\mathbb{P}^2$, $C_1$ is a  conic,%
			\item [$\mathrm{(I.1C})$] $S=\mathbb{P}^2$, $C_1$ is a line,%
			\item [$\mathrm{(I.2.n})$] $S=\mathbb{F}_{n}$ for any $n\ge 0$, $C_1=Z_n$,%
			\item [$\mathrm{(I.3A})$] $S=\mathbb{F}_1$, $C_1\in|2(Z_1+F)|$,%
			\item [$\mathrm{(I.3B})$] $S=\mathbb{F}_1$, $C_1\in|Z_1+F|$,%
			\item [$\mathrm{(I.4A})$] $S=\mathbb{P}^1\times\mathbb{P}^1$,  
			$C_1$ is a $(2,2)$-curve,%
			\item [$\mathrm{(I.4B})$] $S=\mathbb{P}^1\times\mathbb{P}^1$, 
			$C_1$ is a $(2,1)$-curve
			\item [$\mathrm{(I.4C})$] $S=\mathbb{P}^1\times\mathbb{P}^1$, 
			$C_1$ is a $(1,1)$-curve
			\item [$\mathrm{(I.5.m})$] $(S,C)$ is the smooth boundary blow-up of $\mathrm{(I.1A)}$ at $1\le m\le 8$ points  with no 
			three colinear, no six on a conic, and no eight on a cubic smooth away from its single double point with one of the points being that double point,
			\item [$\mathrm{(I.6B.m})$] $(S,C)$ is the smooth boundary blow-up of $\mathrm{(I.1B)}$ at $m\ge 1$ points,%
			\item [$\mathrm{(I.6C.m})$] $(S,C)$ is the smooth boundary blow-up of $\mathrm{(I.1C)}$ at $m\ge 1$ points,%
			\item [$\mathrm{(I.7.n.m})$] $(S,C)$ is the smooth boundary blow-up of $\mathrm{(I.2.n)}$ at $m\ge 1$ points,
			\item [$\mathrm{(I.8B.m})$] $(S,C)$ is the smooth boundary blow-up of $\mathrm{(I.3B)}$ at $m\ge 1$ points,
			\item [$\mathrm{(I.9B.m})$] $(S,C)$ is the smooth boundary blow-up of $\mathrm{(I.4B)}$ at $m\ge 1$ points 
			with no two on the same $(0,1)$-curve,
			\item [$\mathrm{(I.9C.m})$] $(S,C)$ is the smooth boundary blow-up of $\mathrm{(I.4C)}$ at $m\ge 1$ points,
			\item [$\mathrm{(II.1A})$] $S=\mathbb{P}^2$, $C_1$ is a line, $C_2$ is a  conic, 
			\item [$\mathrm{(II.1B})$]  $S=\mathbb{P}^2$, $C_1, C_2$ are  lines,%
			\item [$\mathrm{(II.2A.n})$] $S=\mathbb{F}_n$ for any $n\ge 0$, $C_1=Z_n,\, C_2\in|Z_n+nF|$,
			\item [$\mathrm{(II.2B.n})$]  $S=\mathbb{F}_n$ for any $n\ge 0$, $C_1=Z_n,\, C_2\in|Z_n+(n+1)F|$,%
			\item [$\mathrm{(II.2C.n})$]  $S=\mathbb{F}_n$ for any $n\ge 0$, $C_1=Z_n,\, C_2\in|F|$,%
			\item [$\mathrm{(II.3})$]  $S=\mathbb{F}_1$, 
			$C_1, C_2\in |Z_1+F|$,
			\item [$\mathrm{(II.4A})$]  $S=\mathbb{P}^1\times\mathbb{P}^1$, 
			$C_1,C_2$ are $(1,1)$-curves, 
			\item [$\mathrm{(II.4B})$]  $S=\mathbb{P}^1\times\mathbb{P}^1$,
			$C_1$ is a $(2,1)$-curve, 
			$C_2$ is a $(0,1)$-curve, 
    		\item [$\mathrm{(II.5A.m})$]  $(S,C)$ is
			the smooth boundary blow-up of $\mathrm{(II.1A})$ at $1\le m\le 5$ points
            with no two on the line $\pi(C_1)$, no five on the conic $\pi(C_2)$, and no three collinear,
			
			\item [$\mathrm{(II.5B.m})$]  $(S,C)$ is the smooth boundary blow-up of $\mathrm{(II.1B})$ at $m\ge 1$ points,
			
			\item [$\mathrm{(II.6A.n.m})$] 	$(S,C)$ is the smooth boundary blow-up of $\mathrm{(II.2A.n})$ at $m\ge 1$ points
			with no two on the same curve in $|F|$,
			
			\item [$\mathrm{(II.6B.n.m})$]  $(S,C)$	is the smooth boundary blow-up of $\mathrm{(II.2B.n})$ at $m\ge 1$ points
			with no two on the same curve in $|F|$,
						
			\item [$\mathrm{(II.6C.n.m})$]  $(S,C)$ is the smooth boundary blow-up of $\mathrm{(II.2C.n})$ at $m\ge 1$ points,
			
			\item [$\mathrm{(II.7.m})$]   $(S,C)$	is the smooth boundary blow-up of $\mathrm{(II.3})$ at $m\ge 1$ points
			with no two on the same curve in $|F|$,
					
			\item [$\mathrm{(III.1})$]  $S=\mathbb{P}^2$, $C_1, C_2, C_3$ are lines, 
			
			\item [$\mathrm{(III.2})$]  $S=\mathbb{P}^1\times\mathbb P^1$,
			$C_1, C_2, C_3$ are 
			$(1,1)$-, $(0,1)$-, $(1,0)$-curves, respectively,
			
			\item [$\mathrm{(III.3.n})$]  $S=\mathbb{F}_n$ for any
			$n\ge 0,  \, C_1=Z_n, C_2\in |F|,\, C_3\in |Z_n+nF|$,
			
			\item [$\mathrm{(III.4.m})$]  $(S,C)$ is the smooth boundary blow-up of
			$\mathrm{(III.1})$ at $1\leq m=m_1+m_2+m_3$ points of which $m_i\in\{0,1\}$ lie in $\pi(C_i)$
			for $i\in\{1,2,3\}$, with no three of them collinear,
			
			\item [$\mathrm{(III.5.n.m})$]  $(S,C)$ is the smooth boundary blow-up of
			$\mathrm{(III.3.n})$ at $m\ge 1$ points
			with no two on the same curve in $|F|$,
			and none on the fiber $\pi(C_2)$,
			
			\item [$\mathrm{(IV})$]
			$S=\mathbb{P}^1\times\mathbb{P}^1$, $C_1, C_2$
			are $(1,0)$-curves, $C_3, C_4$ are $(0,1)$-curves.
						
		\end{itemize}
		
	\end{theorem}

\begin{figure}
		\centering
		\includegraphics[width=0.3\textwidth]{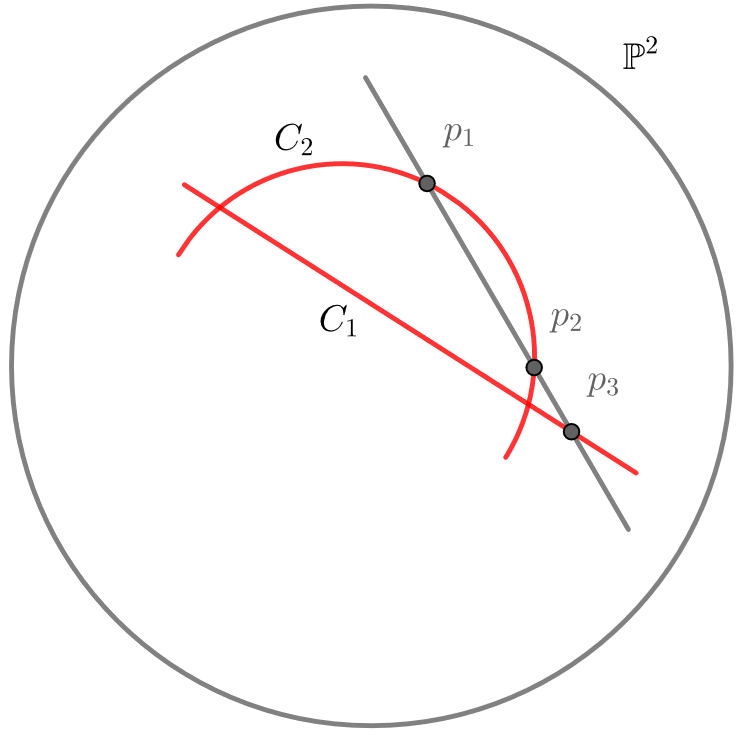}
		\caption{This configuration of points is excluded in (II.5A.3).}
		\label{Figure1}
	\end{figure}

\begin{figure}
		\centering
		\includegraphics[width=0.3\textwidth]{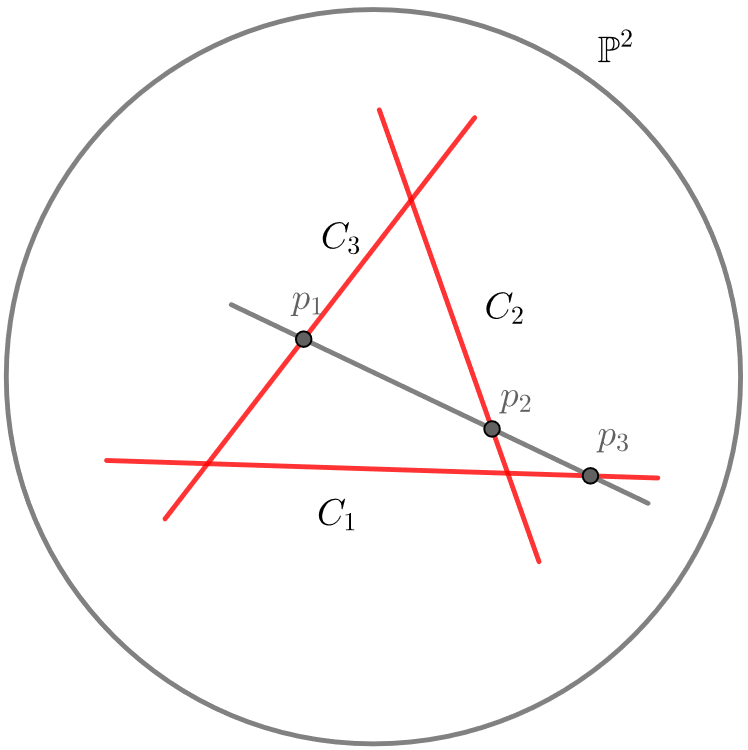}
		\caption{This configuration of points is excluded in (III.4.3).}
		\label{Figure2}
	\end{figure}

\begin{remark} 
{\rm
\label{fixremark}
In Theorem \ref{mainthm} we fixed a few typos and minor omissions from \cite[Theorems 2.1, 3.1]{CR}:

\begin{itemize}
\item In case (I.7.n.m) we fixed the typo (I.2) 
\cite[p. 1253, line 18]{CR} to (I.2.n).
  
\item  Case (II.8.m), $m\in\{1,2,3,4\}$ in \cite[Theorem 3.1]{CR}, while 
correct (modulo the requirement that no two of the points be on the same
(0,1)-fiber\footnote{Note that the generality condition implied by 
SALdP flags with curves not in the boundary (Lemma \ref{ineqwithdeltaALdPcaseLem}) require no two points
on the same fiber and no four on the same $(1,1)$-curve (as 
well as further generality conditions involving more than four points). 
The SALdP
flags with boundary curves (Lemma \ref{bndrycurvesflagslem}) require
no points on the fiber $\pi(C_2)$ and no five points on the 
$(2,1)$-curve $\pi(C_2)$. Thus the latter flags 
already rule out the flag with the $(1,1)$-curve (or any higher bi-degree curves) as such a curve would have to intersect $\pi(C_2)$ but do not rule out the flag
$\{p_1,p_2\}\subset \Sig$ with $\Sig$ a $(0,1)$-fiber different from $\pi(C_2)$.})
coincides with (II.5A.m+1),  $m\in\{1,2,3,4\}$, so we removed the former. 

\item  In case (II.5A.m) we
fixed a typo where $c_1$ and $c_2$ were
interchanged \cite[p. 1260, line $-3$]{CR}.

\item  Still in (II.5A.3), (II.5A.4), and (II.5A.5) 
we added the requirement that no
three of the $m$ points be collinear (such a line would 
intersect the conic $c_1$ at two points and the line $c_2$ at
one point, see Figure \ref{Figure1}). 

\item In (III.4.3) we also added the requirement that the
$m=3$ points not be collinear (such a line would 
intersect each of the lines $c_i$ at one point, see Figure \ref{Figure2}). 

\end{itemize}
} 
\end{remark}

\bremark
\lb{twopointRem}
Since the statement of Theorem \ref{mainthm} stipulates the $C_i$ intersect simply and normally,
the curves in $\mathrm{(II.1A}), \mathrm{(II.4A}), \mathrm{(II.4B}), \mathrm{(II.5A.m})$
intersect at two {\it distinct} points (i.e., we do not allow tangency).
\eremark

The proof of Theorem \ref{mainthm} we give here contains the following main steps.

\med
\noi {\it Step 1.} We classify all \saldp pairs $(s,c)$ with $\rk\Pic(s)\le 2$. 
This is a straightforward computation involving ampleness conditions on 
Hirzebruch surfaces and is carried out in an Appendix (Proposition \ref{prop:classification-rk2-anticanonical}).

\med
\noi {\it Step 2.} We eliminate from the list of Step 1 all pairs that are
obtained from another pair in the list via a {\it smooth boundary blow-up} 
(Definition \ref{properdef}). This is straightforward
and is stated in Proposition \ref{baseprop}.

\med
\noi {\it Step 3.} We introduced the notion of a SALdP flag in
Definition \ref{badcurvebetadef}. In this step (Theorem \ref{MainFlagthm}) we classify
all SALdP flags whose zero-dimensional part lies {\it in the boundary} 
for the pairs in the list of Step 1. This yields the precise
generality conditions appearing in the final classification statement
(Theorem \ref{mainthm}) which slightly improves on the statement in \cite{CR}.

\med
\noi {\it Step 4.} This step was carried out in
Proposition \ref{noawaynotailProp}. 
In sum, we prove, similarly to \cite{CR}, but with different
emphasis, that the list obtained in Step 3 contains {\it all} \saldp pairs. 
The proof of this step is summarized in Figure \ref{Figure0}.

\bpf[Proof of Theorem \ref{mainthm}]
According to \PR{noawaynotailProp} every \saldp pair can be expressed as the 
smooth boundary blow-ups  of a \saldp pairs $(s,c)$ with $\rk\Pic(s)\le 2$.
According to Proposition \ref{baseprop}, Theorem \ref{MainFlagthm},
and Theorem \ref{badflagcriterionprop}, the list in the statement of 
Theorem \ref{mainthm} consists precisely of all \saldp pairs $(s,c)$ with $\rk\Pic(s)\le 2$
that are not the smooth boundary blow-up of any other \saldp pair, together with
all of their smooth boundary blow-ups  that remain \saldpno.

In particular, for the case (I.5.m), there are two types of flags according
to Theorem \ref{MainFlagthm} (I.1A). The second type of flags already
forces $m\le 8$. Feeding this information into the first type
of flags precludes $\Sig$ from being a smooth curve of degree $3$ or higher
and by \cite[Lemma 2.2]{CR} we know the resulting blow-up will be \saldp if
and only if $S$ is del Pezzo.
The del Pezzo--Manin--Hitchin's classification of del Pezzo surfaces \cite{DelPezzo,Manin1,Manin2,Hitchin1975} (see also
\cite{R20} for an exposition in the language of flags close in spirit to this article)
then shows that the only possibilities for $\Sig$ are a smooth line, a smooth conic,
or a singular cubic with a single double point.
The analysis for the cases (II.5A.m) and (III.4.m) is similar but simpler since
the fact that no two points can be on a boundary line component (coming from
the second type of flags in Theorem \ref{MainFlagthm} (II.5A.m), (III.4.m))  
precludes any flag with $\Sig\not\subset c$
of degree 2 or higher. To illustrate this, let us consider the case (II.5A.m).
If $\Sig\not=\pi(C_2)$ is a conic in $\PP^2$ then either it intersects the line $\pi(C_1)$ at two
distinct points which is already taken care of by the first type of flags, or else
it is tangent to $\pi(C_1)$ at a single point, however that is also taken care
of since we cannot blow-up infinitely near points on the smooth part of the boundary
by \cite[Lemma 2.5]{CR}.
Blow-ups of case (I.4A) are contained in the case (I.5.m) by the classification of
del Pezzo surfaces, while blow-ups of (II.4B) are contained in case (II.5A.m) by Remark \ref{fixremark}.
For cases (I.9B.m), (II.6A.n.m), (II.6B.n.m), (II.7.m)
the only obstruction are  pairs of points on the same fiber
according to Theorem \ref{MainFlagthm} (I.4B), (II.2A.n), (II.2B.n), (II.3), and 
for the case (III.5.n.m) one has the previous obstruction as well
as the obstruction of blowing-up any points on the fiber boundary component
of (III.3.n) by Theorem \ref{MainFlagthm} (III.3.n).
For the cases (I.6B.m), (I.6C.m), (I.7.n.m), (I.8B.m), (I.9C.m), (II.5B.m) any smooth boundary blow-ups are 
allowed as by Theorem \ref{MainFlagthm} there are no obstructing SALdP flags
for (I.1B), (I.1C), (I.2.n), (I.3B), (I.4C), (II.1B). 
\epf

\section{Appendix: Classification with small Picard group}

Let $\Fn$ be the $n$-th Hirzebruch surfaces, i.e., the unique
rational surface with Picard group of rank 2 and a curve of self-intersection
$-n$, denoted $Z_n\subset \Fn$. We refer the reader to \cite[\S1.5]{CR}
for our conventions and further background.
Denote by $F$ the class of fiber, i.e., an irreducible smooth rational curve such that $F^2=0$ and $F.Z_n=1$. 
	If $n=0$ when we refer to $Z_0$ and $F$ we
	intend that each is a fiber of a different projection to $\PP^1$.
	Hirzebruch surfaces are ruled toric surfaces  
	and applying adjunction yields 
	\beq
	\label{KFnEq}
	-K_{S}\sim 2Z_n+(n+2)F.
	\eeq
	Recall that every smooth irreducible curve in $|Z_n+nF|$
	(a `zero section') intersects each fiber transversally at a single point
	and does not intersect the `infinity section' $Z_n$.
	Any curve $C$ on $\FF_n$ satisfies 
	\beq\label{Fn-curves}
	C\sim aZ_n+bF,
	\eeq
	with $a,b\in\NN\cup\{0\}$. 
	Also,
	\beq\label{ampleFn}
	\h{$C$ is ample if and
		only if $a>0$ and $b>na$,}
	\eeq
	and furthermore,
	\beq\label{irreducibleFn}
	\h{$C$ is an
		irreducible curve
		only if $C=Z_n$ or $b\ge na\ge0$,}
	\eeq
	and under such conditions the class \eqref{Fn-curves} always contains an irreducible curve which in the latter case is nef.

	\bprop
	\label{prop:classification-rk2-anticanonical}
	Let $S$ be a smooth surface with $\rk(\Pic(S))\le 2$, and let
	$C_1,\ldots,C_r$ be distinct irreducible smooth curves on $S$ such that
	$C=\sum_{i=1}^{r}C_i$  is a divisor with simple normal crossings.
	Then $(S,C)$ is a strongly 
	asymptotically log del Pezzo pair if and only if it
	is one of the following:
	\begin{itemize}
	\item [$\mathrm{(I.1A})$] $S=\mathbb{P}^2$,  $C_1$ is a cubic,%
			\item [$\mathrm{(I.1B})$] $S=\mathbb{P}^2$, $C_1$ is a conic,%
			\item [$\mathrm{(I.1C})$] $S=\mathbb{P}^2$, $C_1$ is a line,%
			\item [$\mathrm{(I.2.n})$] $S=\mathbb{F}_{n}$ for any $n\ge 0$, 
			$C_1=Z_n$,%
			\item [$\mathrm{(I.3A})$] $S=\mathbb{F}_1$, $C_1\in|2(Z_1+F)|$,%
			\item [$\mathrm{(I.3B})$] $S=\mathbb{F}_1$, $C_1\in|Z_1+F|$,%
			\item [$\mathrm{(I.4A})$] $S=\mathbb{P}^1\times\mathbb{P}^1$, 
			$C_1$ is a $(2,2)$-curve,%
			\item [$\mathrm{(I.4B})$] $S=\mathbb{P}^1\times\mathbb{P}^1$, 
			$C_1$ is a $(2,1)$-curve,%
			\item [$\mathrm{(I.4C})$] $S=\mathbb{P}^1\times\mathbb{P}^1$, 
			$C_1$ is a $(1,1)$-curve,%
			\item [$\mathrm{(I.5.1})$] $S=\mathbb{F}_1$, $C_1\in|2Z_1+3F|,$
		\item [$\mathrm{(I.6B.1})$] $S=\mathbb{F}_1$, $C_1\in|Z_1+2F|,$
		\item [$\mathrm{(I.6C.1})$]  $S=\mathbb{F}_1$, $C_1\in|F|,$
		\item [$\mathrm{(II.1A})$]  $S=\mathbb{P}^2$, $C_1$ is a conic, 
		$C_2$ is a line, 
		
		\item [$\mathrm{(II.1B})$]   $S=\mathbb{P}^2$, $C_1,C_2$ are lines,%
		
		\item [$\mathrm{(II.2A.n})$]  $S=\mathbb{F}_n$ for any $n\ge 0$, $C_1=Z_n,\, C_2\in |Z_n+nF|$, 
		
		\item [$\mathrm{(II.2B.n})$] $S=\mathbb{F}_n$ for any $n\ge 0$, $C_1=Z_n,\, C_2\in|Z_n+(n+1)F|$,%
		
		\item [$\mathrm{(II.2C.n})$]  $S=\mathbb{F}_n$ for any $n\ge 0$, $C_1=Z_n,\, C_2\in|F|$,%
		
		\item [$\mathrm{(II.3})$]  $S=\mathbb{F}_1$, $C_1, C_2
		\in|Z_1+F|$,
		
		\item [$\mathrm{(II.4A})$] 
		$S=\mathbb{P}^1\times\mathbb{P}^1$, $C_1,C_2$
		are  $(1,1)$-curves, 
		
		\item [$\mathrm{(II.4B})$] $S=\mathbb{P}^1\times\mathbb{P}^1$, 
		$C_1$ is a $(2,1)$-curve, $C_2$ is a
		$(0,1)$-curve, 
		
		\item [$\mathrm{(II.5A.1 (a)})$]  $S=\mathbb{F}_1$, 
		$C_1\in|2Z_1+2F|,\, C_2\in|F|$, 

        \item [$\mathrm{(II.5A.1 (b)})$]  $S=\mathbb{F}_1$, 
        $C_1\in |Z_1+2F|,\, C_2\in |Z_1+F|$, 
 		
		\item [$\mathrm{(II.5B.1})$]   $S=\mathbb{F}_1$, $C_1\in|F|,
		\, C_2\in |Z_1+F|$,
		
		\item [$\mathrm{(III.1})$]  $S=\mathbb{P}^2$, 
		$C_1, C_2, C_3$ are lines, 
		
		\item [$\mathrm{(III.2})$]  $S=\mathbb{P}^1\times\mathbb P^1$,
		$C_1, C_2, C_3$ are $(1,1)$-, $(0,1)$- and $(1,0)$-curves, respectively,
		
		\item [$\mathrm{(III.3.n})$]  $S=\mathbb{F}_n$ for any
		$n\ge 0,  \, C_1=Z_n, C_2\in |F|,\, C_3\in|Z_n+nF|$,
		
		\item [$\mathrm{(III.4.1})$]  $S=\mathbb{F}_1$, $C_1\in |F|,
 C_2, C_3$ are curves in $|Z_1+F|$,

		\item [$\mathrm{(IV})$]
		$S=\mathbb{P}^1\times\mathbb{P}^1$, $C_1, C_2$
		are  $(1,0)$-curves,  $C_3,C_4$ are  
		$(0,1)$-curves.
		
	\end{itemize}
	
	\eprop

\bremark
The curves in $\mathrm{(II.1B}), \mathrm{(II.4A}), \mathrm{(II.4B}), \mathrm{(II.5A.1) (a)}, \mathrm{(II.5A.1) (b)}$
intersect at two {\it distinct} points, see Remark \ref{twopointRem}.
\eremark
	
	\begin{proof}
		By Castelnuovo's rationality criterion, $S$ is rational, hence by the classification of rational surfaces with $\rk(\Pic(S))\le 2$ it is either $\PP^2$ or $\Fn$ \cite[\S2--3]{CR}. When $S=\PP^2$, $\rk(\Pic(S))=1$, $-K_S\sim 3H$ and we see the possibilities are $\mathrm{(I.1A}),\mathrm{(I.1B}),\mathrm{(I.1C}),\mathrm{(II.1A}),\mathrm{(II.1B}),\mathrm{(III.1})$.
		Assume now that $S=\Fn$. Denote
		$$
		C_i\in|a_iZ_n+b_iF|.
		$$
		Since $-K_S-C$ is nef and using \er{irreducibleFn}, we see that 
		$$
		\sum_i a_i\in\{0,1,2\}, \q 
		\sum_i b_i\in\{0,\ldots,n+2\}. 
		$$
		Note that by \er{irreducibleFn} either $b_i\ge na_i>0$ or else $(a_i,b_i)=(0,1)$, i.e., $C_i$ is a fiber.
		Thus, at most two components of $C_i$ are {\it not } fibers. Also, if $a_i=2$ for some $i$ 
		then $n\le 2$ by $n+2\ge b_i\ge 2n$. So we get the following possibilities when
		$\max_ia_i=2$:
		\beq
		\lb{eqcases1}
		\baeq
		{}[n,(a_1,b_1),\ldots,(a_r,b_r)]\in\big\{
		&[2,(2,4)],[1,(2,3)],[1,(2,2),(0,1)],
		\cr
		&[1,(2,2)],[0,(2,2)],[0,(2,1),(0,1)],
		[0,(2,1)]\big\}.
		\eaeq
		\eeq
		When $\max_ia_i=1$, we split into two subcases:
		when
		there are at least two pairs $(a_i,b_i)$ with all coefficients
		positive:
		\beq
		\lb{eqcases2}
		\baeq
		{}[n,(a_1,b_1),\ldots,(a_r,b_r)]\in\big\{
		&[2,(1,2),(1,2)],[1,(1,2),(1,1)],
		\cr
		&[1,(1,1),(1,1),(0,1)],
		[1,(1,1),(1,1)),(0,(1,1),(1,1)]
		\big\},
		\eaeq
		\eeq
		and otherwise, still with $\max_ia_i=1$, now for all $n$ (so we omit the first index), and
		with $\max_ib_i=n$,
		\beq
		\lb{eqcases3}
		\baeq
		{}[(a_1,b_1),\ldots,(a_r,b_r)]\in\big\{
		&[(1,n),(1,0),(0,1),(0,1)],
		[(1,n),(1,0),(0,1)],
		\cr
		&[(1,n),(1,0)],
		[(1,n),(0,1),(0,1)],
		[(1,n),(0,1)],
		[(1,n)]
		\big\},
		\eaeq
		\eeq
		with $\max_ib_i=n+1$,
		\beq
		\lb{eqcases4}
		\baeq
		{}[(a_1,b_1),\ldots,(a_r,b_r)]\in\big\{
		&[(1,n+1),(1,0),(0,1)],
		\cr
		&[(1,n+1),(1,0)],
		[(1,n+1),(0,1)],
		[(1,n+1)]
		\big\},
		\eaeq
		\eeq
		with $\max_ib_i=n+2$,
		\beq
		\lb{eqcases5}
		\baeq
		{}[(a_1,b_1),\ldots,(a_r,b_r)]\in\big\{
		&[(1,n+2),(1,0)],
		[(1,n+2)]
		\big\},
		\eaeq
		\eeq
		when $\max_ib_i=1$,
		\beq
		\lb{eqcases6}
		\baeq
		{}[(a_1,b_1),\ldots,(a_r,b_r)]
		\in\big\{
		&[(1,0)],
		[(1,0),(0,1)],\ldots,
		[(1,0),\overbrace{(0,1),\ldots,(0,1)}^{n+2}],
		\big\}.
		\eaeq
		\eeq
		and finally when $\max_ia_i=0$,
		\beq
		\lb{eqcases7}
		\baeq
		{}[(a_1,b_1),\ldots,(a_r,b_r)]
		\in\big\{
		&[(0,1)],\ldots,
		[\overbrace{(0,1),\ldots,(0,1)}^{n+2}]
		\big\}.
		\eaeq
		\eeq
		A few of these cases can be eliminated, though most of them
		actually occur. 
		In 
		\er{eqcases1}, [2,(2,4)]
		corresponds to a smooth anticanonical curve in $\FF_2$, 
		which is excluded as $\FF_2$ is not del Pezzo.
		The remaining cases are: $[1,(2,3)]=$  (I.5.1),
		$[1,(2,2),(0,1)]=$  (II.5A.1) (a), 
		$[1,(2,2)]=$  (I.3A),
		$[0,(2,2)]=$  (I.4A),
		$[0,(2,1),(0,1)]=$  (II.4B),
		$[0,(2,1)]=$  (I.4B).
		
		In \er{eqcases2}
		[2,(1,2),(1,2)] is excluded as $Z_2.(Z_2+2F)=Z_2.(2Z_2+4F)=Z_2.(-K_{\FF_2})=0$.
		The remaining cases are: 
		$[1,(1,2),(1,1)]=$  (II.5A.1) (b);\hfill\break  
		$[1,(1,1),(1,1),(0,1)]=$  (III.4.1),
		$[1,(1,1),(1,1))=$  (II.3),
		$(0,(1,1),(1,1)]=$  (II.4A).
		
		In \er{eqcases3},
		$[(1,n),(1,0),(0,1),(0,1)]$
		gives $-K_{\FF_n}-(1-\be_1)(Z_n+nF)
		-(1-\be_2)Z_n
		-(1-\be_3)F-(1-\be_4)F\sim(\be_1+\be_2)Z_n+(n\be_1+\be_3+\be_4)F$,
		that is ample if and only if 
		$n\be_1+\be_3+\be_4> n\be_1+n\be_2$ forcing $n=0$ 
		and this is (IV);
		$[(1,n),(1,0),(0,1)]=$  (III.3.n),
		$[(1,n),(1,0)]=$  (II.2A.n).
		For [(1,n),(0,1),(0,1)] consider 
		$-K_{\FF_n}-(1-\be_1)(Z_n+nF)
		-(1-\be_2)F-(1-\be_3)F\sim(1+\be_1)Z_n+(n\be_1+\be_2+\be_3)F$,
		that is ample if and only if 
		$n\be_1+\be_2+\be_3> n+n\be_1$, i.e., $n=0$, and this is
		(III.3.0).
		For [(1,n),(0,1)], consider 
		$-K_{\FF_n}-(1-\be_1)(Z_n+nF)
		-(1-\be_2)F\sim(1+\be_1)Z_n+(1+n\be_1+\be_2)F$,
		that is ample if and only if 
		$1+n\be_1+\be_2> n+n\be_1$, i.e., $n=0,1$, and these are
		(II.2C.0), (II.5B.1).
		For [(1,n)], $-K_{\FF_n}-(1-\be_1)(Z_n+nF)\sim
		(1+\be_1)Z_n+(2+n\be_1)F$ implying $n=0,1,2$ and
		these are (I.2.0), (I.3B), while the case $n=2$ is excluded
		as in the first paragraph.
		
		In \er{eqcases4}:
		[(1,n+1),(1,0),(0,1)] 
		$-K_{\FF_n}-(1-\be_1)(Z_n+(n+1)F)
		-(1-\be_2)Z_n-(1-\be_3)F\sim(\be_1+\be_2)Z_n+((n+1)\be_1+\be_3)F$
		that is ample if and only if
		$(n+1)\be_1+\be_3>n(\be_1+\be_2)$, forcing $n=0$,
		and this is (III.2);
		$[(1,n+1),(1,0)]=$  (II.2B.n);
		[(1,n+1),(0,1)],
		$-K_{\FF_n}-(1-\be_1)(Z_n+(n+1)F)
		-(1-\be_2)F\sim(1+\be_1)Z_n+((n+1)\be_1+\be_2)F$
		that is ample if and only if
		$(n+1)\be_1+\be_2>n(1+\be_1)$, i.e., $n=0$
		and this is (II.2B.0);
		[(1,n+1)],
		$-K_{\FF_n}-(1-\be_1)(Z_n+(n+1)F)
		\sim(1+\be_1)Z_n+(1+(n+1)\be_1)F$
		that is ample if and only if
		$1+(n+1)\be_1>n(1+\be_1)$, i.e., $n=0,1$,
		and these are (I.4C), (I.6B.1).
		
		In \er{eqcases5},
		[(1,n+2),(1,0)],
		$-K_{\FF_n}-(1-\be_1)(Z_n+(n+2)F)
		-(1-\be_2)Z_n\sim(\be_1+\be_2)Z_n+(n+2)\be_1F$,
		i.e., $2\be_1>n\be_2$, so $n=0$ and this is (II.4B);
		[(1,n+2)],
		$-K_{\FF_n}-(1-\be_1)(Z_n+(n+2)F)
		\sim(1+\be_1)Z_n+(n+2)\be_1F$,
		i.e., $(n+2)\be_1>n+n\be_1$, i.e., $n=0$ and
		this is (I.4B).
		
		In \er{eqcases6},
		there is one $Z_n$ and $k$ fibers, with $k\ge0$:
		$-K_{\FF_n}-(1-\be_1)Z_n
		-(1-\be_2)F
		-\ldots
		(1-\be_{k+1})F
		\sim(1+\be_1)Z_n+(n+2-k+\be_2+\ldots+\be_{k+1})F$,
		that is ample if and only if 
		$
		n(1+\be_1)
		<
		n+2-k+\be_2+\ldots+\be_{k+1}
		$,
		i.e., $k=0,1,2$.
		When $k=0$ this is (I.2.n), when $k=1$ this is
		(II.2C.n), and when $k=2$ this 
		means $n\be_1<\be_2+\be_3$, forcing $n=0$ and this
		is (III.3.0). 
		
		Finally, in \er{eqcases7},
		there are $k>0$ fibers, so
		$-K_{\FF_n}-(1-\be_1)F
		-\ldots
		(1-\be_{k})F
		\sim2Z_n+(n+2-k+\be_1+\ldots+\be_{k})F$,
		that is ample if and only if 
		$
		2n
		<
		n+2-k+\be_1+\ldots+\be_{k}
		$,
		i.e., $n=k=1$ and (I.6C.1), or $n=0$ and $k=1,2$
		and these are (I.2.0), (II.2A.0).
		\epf

As a corollary, we obtain the following which is needed as a step in the proof
of Theorem \ref{mainthm}.

\bprop
\lb{baseprop}
Let $\big(S,C=\sum_{i=1}^{r}C_i\big)$ be \saldp pair with $\rk(\Pic(S))\le 2$ that is not the smooth boundary blow-up (Definition \ref{properdef}) of any other \saldp pair.
	Then $(S,C)$ is a strongly 
	asymptotically log del Pezzo pair if and only if it
	is one of the following:
	\begin{itemize}
	\item [$\mathrm{(I.1A})$] $S=\mathbb{P}^2$,  $C_1$ is a cubic,%
			\item [$\mathrm{(I.1B})$] $S=\mathbb{P}^2$, $C_1$ is a conic,%
			\item [$\mathrm{(I.1C})$] $S=\mathbb{P}^2$, $C_1$ is a line,%
			\item [$\mathrm{(I.2.n})$] $S=\mathbb{F}_{n}$ for any $n\ge 0$, 
			$C_1=Z_n$,%
			\item [$\mathrm{(I.3A})$] $S=\mathbb{F}_1$, $C_1\in|2(Z_1+F)|$,%
			\item [$\mathrm{(I.3B})$] $S=\mathbb{F}_1$, $C_1\in|Z_1+F|$,%
			\item [$\mathrm{(I.4A})$] $S=\mathbb{P}^1\times\mathbb{P}^1$, 
			$C_1$ is a $(2,2)$-curve,%
			\item [$\mathrm{(I.4B})$] $S=\mathbb{P}^1\times\mathbb{P}^1$, 
			$C_1$ is a $(2,1)$-curve,%
			\item [$\mathrm{(I.4C})$] $S=\mathbb{P}^1\times\mathbb{P}^1$, 
			$C_1$ is a $(1,1)$-curve,%
		\item [$\mathrm{(II.1A})$]  $S=\mathbb{P}^2$, $C_1$ is a conic, 
		$C_2$ is a line, 
		
		\item [$\mathrm{(II.1B})$]   $S=\mathbb{P}^2$, $C_1,C_2$ are
		lines,%
		
		\item [$\mathrm{(II.2A.n})$]  $S=\mathbb{F}_n$ for any $n\ge 0$, $C_1=Z_n,\, C_2\in |Z_n+nF|$, 
		
		\item [$\mathrm{(II.2B.n})$] $S=\mathbb{F}_n$ for any $n\ge 0$, $C_1=Z_n,\, C_2\in|Z_n+(n+1)F|$,%
		
		\item [$\mathrm{(II.2C.n})$]  $S=\mathbb{F}_n$ for any $n\ge 0$, $C_1=Z_n,\, C_2\in|F|$,%
		
		\item [$\mathrm{(II.3})$]  $S=\mathbb{F}_1$, $C_1, C_2\in
		|Z_1+F|$,
		
		\item [$\mathrm{(II.4A})$] 
		$S=\mathbb{P}^1\times\mathbb{P}^1$, $C_1,C_2$
		are  $(1,1)$-curves, 
		
		\item [$\mathrm{(II.4B})$] $S=\mathbb{P}^1\times\mathbb{P}^1$, 
		$C_1$ is a $(2,1)$-curve, $C_2$ is a
		$(0,1)$-curve, 

		\item [$\mathrm{(III.1})$]  $S=\mathbb{P}^2$, 
		$C_1, C_2, C_3$ are lines, 
		
		\item [$\mathrm{(III.2})$]  $S=\mathbb{P}^1\times\mathbb P^1$,
		$C_1$ is a $(1,1)$-curve, $C_2$ is a $(0,1)$-curve, $C_3$ is a $(1,0)$-curve, 
		
		\item [$\mathrm{(III.3.n})$]  $S=\mathbb{F}_n$ for any
		$n\ge 0,  \, C_1=Z_n, C_2\in |F|,\, C_3\in|Z_n+nF|$,	
		
		\item [$\mathrm{(IV})$]
		$S=\mathbb{P}^1\times\mathbb{P}^1$, $C_1, C_2$
		are $(1,0)$-curves,  $C_3,C_4$ are  
		$(0,1)$-curves,
		
	\end{itemize}

\eprop

\bpf 
By \PR{prop:classification-rk2-anticanonical}, $(s,c)$ must
be one of the pairs listed there. Of those listed there, 
(I.5.1), (I.6B.1), (I.6C.1), (II.5A.1 (a)), (II.5A.1 (b)), (II.5B.1), (III.4.1)
are manifestly obtained as the 
smooth boundary		 
blow-up (Definition \ref{properdef}) of the pairs (I.1A), (I.1B), (I.1C),  (II.1A) (blow-up on the line in (a), blow-up on the conic
in (b)), (II.1B),
(III.1),
respectively. This completes the proof since the list of \PR{baseprop} coincides
with that of \PR{prop:classification-rk2-anticanonical} modulo those 7 cases.
\epf

\begin{spacing}{0}

\def\bi{\bibitem}

	\end{spacing}
	
	\bigskip	
	
	{\sc University of Maryland }
	
	{\tt yanir@alum.mit.edu}

\end{document}